\documentclass[3p,times]{elsarticle}
\usepackage{amssymb}
\usepackage{amsmath} 
\usepackage{amsthm}
\usepackage{soul,xcolor}
\usepackage{caption}
\usepackage{subcaption}
\usepackage{graphicx}
\usepackage{tikz}
\usepackage{pgfplots}
\usetikzlibrary{arrows.meta, patterns, shapes, calc}

\usepackage{pgfplots}
\usepackage{color}
\usepackage[hidelinks]{hyperref}
\usepackage{multirow}
\usepackage{upgreek}
\usetikzlibrary{arrows.meta}

\newcommand{\vertiii}[1]{{\left\vert\kern-0.1ex\left\vert\kern-0.1ex\left\vert #1 \right\vert\kern-0.1ex\right\vert\kern-0.1ex\right\vert}} 

\makeatletter
\def\hlinewd#1{%
  \noalign{\ifnum0=`}\fi\hrule \@height #1 \futurelet
   \reserved@a\@xhline}
\makeatother

\allowdisplaybreaks

\def\el{{\nonumber}}

\usepackage{regexpatch} 
\def\ferrorp{{\varepsilon'_{\rm p}}}
\def\ferrord{{\varepsilon'_{\rm d}}}
\def\mM{\hbox{{\sffamily M}}}
\def\mK{\hbox{{\sffamily K}}}
\def\mP{\hbox{{\sffamily P}}}
\def\mQ{\hbox{{\sffamily Q}}}
\def\mF{\hbox{{\sffamily F}}}
\def\mD{\hbox{{\sffamily D}}}
\def\mz{\hbox{{\sffamily z}}}
\def\matu{\hbox{{\sffamily u}}}

\def\matxi{\upxi}
\def\cl{\nonumber\\}
\def\el{\nonumber}

\begin{document}

\begin{frontmatter}

\title{A variational multiscale approach to goal-oriented\\ error estimation in finite element analysis of \\ convection-diffusion-reaction equation problems}

\author[deca,cimne]{Sheraz Ahmed Khan}
\ead{sheraz.ahmed@upc.edu}

\author[deca,cimne]{Ramon Codina\corref{ca}}
\ead{ramon.codina@upc.edu}

\author[ovgu]{Hauke Gravenkamp}
\ead{hauke.gravenkamp@ovgu.de}

\address[deca]{Universitat Polit\`{e}cnica de Catalunya (UPC), Jordi Girona 1-3, Edifici C1, Barcelona, 08034, Spain}
\address[cimne]{International Centre for Numerical Methods in Engineering (CIMNE), Campus Nord UPC, Gran Capit\`a s/n,
Barcelona, 08034, Spain}
\address[ovgu]{Institute of Materials, Technologies and Mechanics, Otto von Guericke University Magdeburg, Germany}

\cortext[ca]{Corresponding author}

\begin{abstract}
This paper presents a goal-oriented a posteriori error estimation framework for linear functionals in the stabilized finite element discretization of the stationary convection-diffusion-reaction (CDR) equation. The theoretical framework for error estimation is based on the variational multiscale (VMS) concept, where the solution is decomposed into resolved (finite element) and unresolved (sub-grid) scales. In this work, we propose an orthogonal sub-grid scale (OSGS) method for a goal-oriented error estimation in VMS discretizations. In the OSGS approach, the space of the sub-grid scales (SGSs) is orthogonal to the finite element space. The error is estimated in the quantity of interest, given by the linear functional $Q(u)$ of the unknown $u$. If the SGS $u'$ is estimated, the error in the quantity of interest can be approximated by $Q(u')$. Our approach is compared with a duality-based a posteriori error estimation method, which requires the solution of an additional auxiliary problem. The results indicate that both methods yield similar error estimates, whereas the VMS-based explicit approach is computationally less expensive than the duality-based implicit approach. Numerical tests demonstrated the effectiveness of our proposed error estimation techniques in terms of the quantity of interest functionals. 
\end{abstract}


\end{frontmatter}

\section{Introduction}

Over the past few decades, finite element (FE) methods have become a cornerstone for solving complex mathematical and engineering problems. However, conventional FE approaches often struggle to accurately capture solutions to boundary value problems, in particular those dominated by convection effects. Reliable numerical results require not only robust computational methods but also careful estimation of the approximation errors. In this work, we focus on quantifying these errors by performing a goal-oriented a posteriori error estimation for a specific quantity of interest, using a stabilized finite element formulation of the stationary convection-diffusion-reaction (CDR) equation.

Stabilized FE formulations often incorporate mechanisms for a posteriori error estimation, enabling adaptive strategies to reduce the discretization errors and enhance both accuracy and computational efficiency in complex simulations. Considering this, we employ stabilization techniques such as the orthogonal sub-grid scale (OSGS) method \cite{codina2000stabilization} to develop the stabilized formulation of the primal and the dual problem. The OSGS approach is derived within the overarching framework of the variational multiscale method (VMS), established by the pioneering work of Hughes et al.~\cite{hughes1995multiscale,hughes1998variational}. These stabilized methods are based on the introduction of additional terms to the standard Galerkin FE method. In general, the VMS method consists of splitting a continuous solution into a resolved scale (FE solution) and an unresolved or sub-grid scale (SGS), which can be understood as the numerical error of the resolved scale. The resolved and unresolved scales are also known as coarse and fine scales, respectively. In the OSGS approach, the latter are taken as orthogonal to the former. The addition of a numerical error or fine-scale part to the resolved scales recovers the part that is not represented by the FE solution. This approach has been extensively used as a foundational framework for deriving stabilization techniques \cite{hughes2018multiscale,codina2018variational}. In recent years, the VMS method has also served as a useful tool for estimating the a posteriori error in fluid mechanics problems \cite{hauke2023review}. 

Independent of the treatment of instabilities, a posteriori error estimation is an active area of research \cite{ainsworth1997posteriori,chamoin2023introductory,gratsch2005posteriori}, and it is increasingly being used as a tool to improve the quality of FE solutions. In this paper, we employ concepts of a posteriori error estimation that are also based on the VMS theory. In fact, the initial development of the VMS method \cite{hughes1998variational} already proposed this concept be utilized as a framework for this purpose. If $u$ is the unknown of the problem and $u_h$ the FE solution, this idea was based on modeling the SGS $u'=u-u_h$ directly as an a posteriori error estimator. Since then, this concept has been explored for the CDR equations by Hauke et al.~\cite{hauke2006multiscale,hauke2006proper}, for the Poisson equations by Larson and M{\aa}lqvist \cite{larson2005adaptive,larson2007adaptive}, and for the compressible Navier-Stokes equations by Bayona-Roa et al.~\cite{bayona2018variational}. In another study, Colom{\'e}s et al.\ investigated the application of the VMS error estimation to error propagation in uncertainty quantification \cite{colomes2018robustness}. The use of the SGSs as an a posteriori error estimator for transport problems was investigated by Hauke et al.~\cite{hauke2008variational}. Another application was explored by Codina et al.\ for reduced-order models \cite{codina2021posteriori}. Based on the above ideas, error estimates for solid mechanics problems were proposed by Baiges et al.~\cite{baiges2017variational}. 

In general, a posteriori error estimation techniques can be classified into two sub-categories: explicit and implicit methods. Explicit methods are based on direct post-processing of FE solutions. This type of error estimators were introduced by Babu{\v{s}}ka and Rheinboldt \cite{babuvska1979analysis}, and later extended to two dimensions \cite{kelly1983posteriori,babuvska1987feedback}. On the other hand, implicit methods originate from the class of error estimates that require solving an additional auxiliary problem to estimate the error. 

In practical applications, researchers are often more concerned with the error computation of a specific quantity of interest rather than the error estimation in a certain norm. The implicit error estimation framework plays a crucial role in designing goal-oriented adaptive strategies, where it enables precise control by incorporating information from both the primal and dual problems. The representative work on the implicit error estimation was initiated by \cite{babuvska1978posteriori,babuvvska1978error}, with significant contributions later presented in \cite{jin1998posteriori,diez1998posteriori}. For Stokes problems, a study by Larsson et al.~\cite{larsson2010flux} presented the implicit error estimates using mixed FE formulations, where they employed an energy-like norm to measure the error. For convection-diffusion problems, the use of an energy-like norm yields formulations that are not robust with respect to convection (see, e.g. \cite{korotov-2008,zhang-et-al-2011} and the discussion in \cite{john-novo-2013,du2021robust,sharma2021robust,tobiska2015robust}). We explain in this paper the difficulties of using implicit methods based on global a posteriori error estimates to derive goal-oriented error estimates.

In this work, we propose using the OSGS method to conduct a goal-oriented a posteriori error estimation for the VMS discretization of the stationary CDR equation. Our objective is to present a framework for the OSGS method that provides highly accurate error predictions for a user-defined quantity of interest. The use of the OSGS formulation and the assessment of its numerical performance are the main contributions of this paper.

The quantity of interest is assumed to be the result of a linear functional $Q(u)$ applied to the solution $u$, and the goal is to estimate the error of this value. Furthermore, we present an implicit and an explicit method for the goal-oriented a posteriori error estimation using the OSGS method. In the explicit approach, once the SGS $u'$ is estimated, the error in the quantity of interest can be approximated as $Q(u')$. In contrast to the explicit approach, the implicit method is based on the duality principle, which requires solving an additional adjoint problem to estimate the error in the quantity of interest functional. Initially, Hauke and Fuster \cite{hauke2009variational} derived explicit error estimates for the quantity of interest functional for transport problems by using the VMS technology. Later, a similar approach was followed by Granzow et al.\ \cite{granzow2017output}, in which they compared the efficiencies of the explicit and the implicit error estimates for the CDR equations. Garg and Stonger \cite{garg2019local} utilized the SGS information to derive adjoint error estimates for non-trivial quantities of interest. Among the different approaches for adjoint-based goal-oriented error estimation techniques, we highlight the works for elliptic problems by Abdulle and Nonnenmacher \cite{abdulle2013posteriori} and Wildey et al.~\cite{wildey2008posteriori}, and for non-linear reaction-diffusion problems by Li and Yi \cite{li2022posteriori}. Several authors have applied a goal-oriented error estimation for the FE discretization of the convection-diffusion problems, as explored in previous studies \cite{kuzmin2010goal,cnossen2006aspects,valseth2020goal,chaudhry2014enhancing}. 

The remainder of the paper is organized as follows. Section 2 introduces the strong and weak formulations of the model problem. In addition, it includes a description of the VMS framework for designing and implementing the stabilized FE method we use. Section 3 states the strong and weak forms of the dual problem along with the stabilization and implementation procedure. In Section 4, the explicit and implicit error estimators are presented, and the equivalence between them is proved. Section 5 discusses numerical tests that evaluate the performance of error estimators for a specific quantity of interest. Finally, Section 6 summarizes the key findings of this study. 

\section{Primal problem}

\subsection{Strong form of the primal problem}

Let $\Omega$ be a bounded domain in $\mathbb{R}^{n_{\text{sd}}}$, where $n_{\text{sd}}$ denotes the number of spatial dimensions of the problem. In this work, we consider one and two-dimensional numerical examples, but the concepts also apply to $n_{\text{sd}} = 3$. 

The strong form of the stationary CDR equation consists of finding $u:\Omega \rightarrow \mathbb{R}$ such that for a given function $f:\Omega \rightarrow \mathbb{R}$, the following equations are satisfied
\begin{align}\label{strprim}
\begin{split}
       -k \Delta u + a \cdot \nabla u + su &= f \qquad \text{in} \quad  \Omega \\
      \hspace{2.5cm} u &= 0 \qquad  \text{on} \quad \Gamma
\end{split}
\end{align}
where $k$ is the diffusion coefficient, $a$ is the convection velocity, $s$ is the reaction coefficient, and $f$ is a body load. We assume that $k > 0$, $s\geq 0$, and $a$ is divergence-free (for simplicity), continuous, and bounded. For the sake of simplicity, we consider only homogeneous Dirichlet conditions on $\Gamma$.

To simplify notation, let us introduce the linear differential operator $\mathcal{L}$ as
\begin{equation} \label{strprimeqnot}
  \mathcal{L}u = -k \Delta u + a \cdot \nabla u + su
\end{equation}

\subsection{Weak form of the primal problem and Galerkin finite element approximation}

In order to present the weak form of the problem, let us introduce the space for the trial functions and the test functions, $V = H^1_0(\Omega)$, i.e.,  
$ V=\left\{ v \in H^{1} (\Omega) ~ : ~ v|_{\Gamma} = 0  \right\}$. 
The weak form is derived by multiplying the test function $v \in V$ by the strong form of the problem and then integrating by parts. Thus, a variational formulation can be constructed as follows: find $u \in V$ such that 
\begin{equation} \label{varform}
	B(u,v)=L(v) \quad \forall v\in V 
\end{equation}
where $B(\cdot,\cdot)$ represents the bilinear form of the problem, and $L$ represents the linear form arising from the forcing term. In our problem, these are given by
\begin{align*}
  B(u,v) &=\int_{\Omega}{k \nabla u\cdot \nabla v} \ \text{d}\Omega + \int_{\Omega}{a \cdot \nabla u v} \ \text{d}\Omega + \int_{\Omega}{s u v} \ \text{d}\Omega\\
  L(v) &= \int_{\Omega}{f v} \ \text{d}\Omega
\end{align*}

Let us consider now an FE partition ${\cal T}_h = \{ K\}$ of the domain $\Omega$, with $h_K = {\rm diam}({K})$ and 
$h = \max \{ h_K : K\in {\cal T}_h\}$. 
For the sake of simplicity, we will consider the FE mesh as quasi-uniform, $h$ being the characteristic element size. 
Let  $V_h \subset V$ be a conforming finite space, constructed as
\begin{align*}
  V_h &=\left\{ v_h \in V \ | \ v_h\vert_K \in \mathbb{P}_p(K) \right\}
\end{align*}
where $\mathbb{P}_p(K)$ represents the space consisting of polynomials of degree $p$ in subdomain $K$. Thus, the standard Galerkin approximation of \eqref{varform} reads: find $u_h \in V_h$ such that 
\begin{equation*} 
	B(u_h,v_h)=L(v_h)\quad \forall v_h\in V_h 
\end{equation*}

\subsection{Variational multiscale formulation of the primal problem}

In the following, we briefly summarize the stabilization approach that has been extensively employed for the CDR equation in previous works \cite{codina2000stabilized,codina1998comparison,codina2000stabilization,principe2010stabilization}. This summary is required, since, as we shall see, there are key steps in developing the formulation that will be used thereafter. A comprehensive overview of this method can be found in \cite{codina2018variational}. The foundation is based on the VMS framework \cite{hughes1995multiscale,hughes1998variational}, which consists of decomposing the continuous space into coarse and fine scales. Thus, we have 
\begin{align}
  V &=V_h \oplus V'
\end{align}
Here, $V_h$ represents the coarse-scale space for the trial and test functions,  whereas $V'$ denotes their corresponding fine-scale space. 
Once the FE spaces are established, we further split the trial functions $u$ and test functions $v$ into coarse and fine scales:
\begin{align} \label{subdecom}
\begin{split}
  u&=u_h+u' \qquad u_h \in V_h, ~ u' \in V'  \\
  v&=v_h+v' \qquad v_h \in V_h, ~ v' \in V' 
\end{split}
\end{align}
where $u, u_h$, and $u'$ represent the exact solution, coarse and fine scales, respectively, and likewise for the test functions.
\par
The decomposition outlined above can be incorporated into the variational formulation \eqref{varform}. Thanks to the bilinearity of $B(\cdot, \cdot)$, the weak form transforms into two subproblems: find $u_h \in V_h$ and $u' \in V'$ such that 
\begin{align}
  B(u_h,v_h)+B(u',v_h) & = L(v_h) \qquad \forall v_h\in V_{h} \label{coarsescale} \\ 
  B(u_h,v')+B(u',v') & = L(v') \qquad \forall v'\in V' \label{finescale}
\end{align}
These two subproblems represent the weak form for the coarse and fine scales, respectively. The second term on the left-hand side in Eq.~\eqref{coarsescale} is the contribution of the SGS $u'$ to the FE component, which can be obtained from the fine-scale problem \eqref{finescale}. This yields the stabilized FE formulation in the VMS framework.

We wish to obtain an {\em approximate} expression for the SGS $u'$. Thus, the first assumption we make is that it vanishes on the element boundaries. Even though this can be relaxed (see \cite{codina2009subscales}), this is enough for the problems in which we are interested. However, in diffusion-dominated problems it is crucial to have an approximation of $u'$ on the element boundaries for a posteriori error estimation~\cite{codina-gravenkamp-khan-2025}. With this assumption, if we integrate the second term on the left-hand side of Eq.~\eqref{coarsescale} by parts, we get:
\begin{equation}
  \label{coarsescale2}
  B(u_h,v_h)+\sum_K{{\left\langle u',{\mathcal{L}}^*v_h\right\rangle }_K}=L(v_h)
\end{equation}
for all $v_h \in V_h$. Here ${{\left\langle \cdot, \cdot \right\rangle }_\omega}$ indicates the integral of the product of two functions over a subdomain $\omega\subset\Omega$ (the subscript will be omitted when $\omega = \Omega$), and $\mathcal{L}^*$ denotes the adjoint of the differential operator $\mathcal{L}$. This adjoint operator $\mathcal{L}^*$ arises as a result of applying integration by parts to the differential operator $\mathcal{L}$. Recalling that we have taken $\nabla\cdot a = 0$, and that $u'$ is considered negligible on the element boundaries, $\mathcal{L}^*$ is given by 
\begin{equation}
  {\mathcal{L}}^*v =-k \Delta v - a \cdot \nabla v+sv \label{adj-op}
\end{equation}

Among the family of VMS methods, various approaches exist with the key variation being in the approximation of the fine-scale problem \eqref{finescale}. At this point, additional simplifications are necessary to derive a computationally efficient numerical method. To obtain the SGS approximation in terms of the FE component, and taking also the test function $v'$ as zero on the element boundaries, the fine-scale problem \eqref{finescale} can be rewritten as
\begin{align}
  B(u',v') &= L(v')-B(u_h,v')\nonumber\\
           &= L(v')-\sum_K{{\left\langle {\mathcal{L}u_h},v'\right\rangle }_K}\nonumber\\
           &= \sum_K{{\left\langle {\mathcal{R}u_h},v'\right\rangle }_K} \nonumber\\
           & = \sum_K \langle {\cal L}u',v'\rangle_K\nonumber
\end{align}
where $\mathcal{R}u_h=f-\mathcal{L}u_h$ is the FE residual. If $P'$ denotes the (broken)\footnote{A \textit{broken} projection, in this context, means that it is defined element-wise, without enforcing regularity across element boundaries.} $L^2$ projection onto $V'$, this equation states that
$$ P' ({\cal L}u') = P'(\mathcal{R}u_h)$$

At this point, the critical approximation is the introduction of the stabilization parameter $\tau$ defined element-wise, such that $\tau^{-1}u'$ approximates, in a certain sense, ${\cal L}u$ (see \cite{codina9} for a motivation based on a Fourier analysis). Then
$$ P'(\mathcal{R}u_h) = P' ({\cal L}u') \approx P' (\tau^{-1}u') =  \tau^{-1}u' $$
Therefore, the approximation we propose is
\begin{equation}
  u'|_K \approx u'_h\vert_K := \tau P'(\mathcal{R}u_h)|_K \label{SGSappr}
\end{equation}
where the approximate SGS $u'_h$ can be computed element by element, and the stabilization parameter $\tau$  is given by \cite{codina2000stabilization}:
\begin{align} 
  \tau =\left( c_1\frac{k}{h^{2}}+c_2\frac{|a|}{h}+c_3s \right)^{-1} \label{tauelem}
\end{align}
where $c_1, c_2$ and $c_3$ are algorithmic constants that, for linear elements, we take as
\begin{equation*}
  c_1=4, \ c_2=2, \ c_3=1
\end{equation*}
In \eqref{tauelem}, $h$ is the diameter of the element under consideration and $\vert a \vert$ a characteristic velocity in the element (for example $\vert a \vert = \Vert a \Vert_{L^\infty(K)}$). However, to lighten the notation, we will omit the subscript referring to the element, as it has been done for $\tau$ in \eqref{tauelem}.

Following \eqref{SGSappr}, this approximation allows us to rewrite the coarse-scale problem \eqref{coarsescale2} by modifying the SGSs component in terms of the FE variables. This yields the stabilized formulation: find $u_h\in V_h$ such that
\begin{equation}
  B(u_h,v_h)+\sum_K{{\left\langle \tau P'(\mathcal{R}u_h),{\mathcal{L}}^*v_h\right\rangle }_K}=L(v_h)\quad \forall v_h\in V_h
\end{equation}
The introduced stabilizing term can be written as
\begin{align}
  \sum_K{{\left\langle \tau P'(\mathcal{R}u_h),{\mathcal{L}}^*v_h\right\rangle }_K}
  & = L(v_h) - B(u_h,v_h) \nonumber \\
  & = \sum_K\left\langle \mathcal{R}u_h,v_h\right\rangle_K - \sum_K \left\langle k \partial_n u_h,v_h\right\rangle_{\partial K}
   \label{eqnotstabform}
\end{align}
where $\partial_n$ stands for derivative in the direction normal to $\partial K$. In this expression, derivatives are understood in the classical sense. However, we could also consider second derivatives in the sense of distributions, and simply write
\begin{align}
L(v_h) - B(u_h,v_h) = \left\langle \mathcal{R}u_h,v_h\right\rangle  \label{eqnotstabform-bis}
\end{align}
where now $\left\langle \cdot ,\cdot \right\rangle$ is used for the duality pairing based on the integral between $V$ and its dual, where  $\mathcal{R}u_h$ belongs. Note that, using this notation,
\begin{align*}
B(u,v) = \langle {\cal L} u , v \rangle = \langle u , {\cal L}^{*} v \rangle
\end{align*}
for all $u, v\in V$.

At this point, it only remains to choose $P'$. A simple choice is to take $P'=I$. When acting on the FE residual, this results in a variant known as the Algebraic SGS (ASGS) method, used for example in \cite{hughes1995multiscale}. However, the approach we follow is based on selecting the SGS space in a specific fashion such that it is orthogonal to the FE space, i.e., we consider the decomposition $V=V_h + V^{\bot}_h$. This implies $P'=P^\bot=I-P_h$, where $P_h$ denotes the broken $L^2$ projection onto the FE space $V_h$. This choice of the SGS leads to the orthogonal SGS (OSGS) method \cite{codina2000stabilization}. Thus,  the representation of the SGS $u'$ for the OSGS method can be written as:
\begin{equation} \label{subgridOSGS}
  u' \approx u'_h:= \tau P^{\bot}(\mathcal{R}u_h)\quad \hbox{(computed element-wise)}
\end{equation}

The final stabilized formulation for the OSGS method turns into: find $u_h \in V_h$ such that
\begin{equation}\label{stabformprim}
  B(u_h,v_h)+\sum_K{{\left\langle \tau P^{\bot}(\mathcal{R}u_h),{\mathcal{L}}^*v_h\right\rangle }_K}=L(v_h)\quad \forall v_h \in V_h
\end{equation}
Expression \eqref{subgridOSGS} is the approximation we propose for the SGS, and we have distinguished it from the exact SGS (arising from a given splitting $V = V_h\oplus V'$) by the subscript $h$. We shall see that this distinction is necessary for the following developments. If we take the discrete solution to be $u_h + u_h'$, there will be an error with respect to the exact solution $u$ that we call $\ferrorp$, i.e.,
\begin{align}
u = u_h + u'_h+ \ferrorp \label{ferrorp}
\end{align}
Obviously, the FE function $u_h$ appearing in this expression is not the same as in \eqref{subdecom}, but distinguishing between them will not be necessary. The subscript in $\ferrorp$ refers to the fact that we are considering the primal problem.

\subsection{Implementation of the OSGS formulation for the primal problem}

In this subsection, we describe the implementation procedure for the OSGS-based stabilized method in a matrix version. Recalling the concept of the OSGS method, the SGS $u'_h= \tau P^{\bot}(f-\mathcal{L}u_h)$ is computed as the residual's component orthogonal to the FE space, i.e., $P^{\bot}=I-P_h$. Here, we denote the broken projection $P_h$ of the FE residual as
\begin{equation*}
  \xi_h = P_h(f-\mathcal{L}u_h)
\end{equation*}
The computation of $\xi_h$ is achieved by solving
\begin{equation} \label{projcomp}
  \sum_K \langle \mathcal{L}u_h, \zeta_h\rangle_K + \sum_{K} {\left\langle \xi_h,\zeta_h\right\rangle}_K =\left\langle f,\zeta_h\right\rangle
\end{equation}
for all $\zeta_h \in V_h$. Combining Eqs.~\eqref{stabformprim} and \eqref{projcomp}, the stabilized formulation can be expressed as a system of equations given by:
\begin{align}
  B(u_h,v_h)-\sum_K \tau {\left\langle \mathcal{L}u_h, \mathcal{L}^{*}v_h \right\rangle}_K - \sum_K \tau {\left\langle \xi_h, \mathcal{L}^{*}v_h \right\rangle}_K &= L(v_h) - \sum_{K}\tau {\left\langle f,\mathcal{L}^{*}v_h\right\rangle}_K \label{disceq1}\\
  \sum_K \langle \mathcal{L}u_h, \zeta_h\rangle_K + \sum_{K} {\left\langle \xi_h,\zeta_h\right\rangle}_K &= L(\zeta_h)\label{disceq2}
\end{align}
In this study, the projections $\xi_h$ are computed implicitly following the method outlined in \cite{codina2008analysis}. Let us rewrite the discrete form of the coupled system of equations \eqref{disceq1} and \eqref{disceq2}, leading to the matrix representation
\[
\begin{bmatrix}
\mK & -\mP_{\tau} \\
\mD & \mM
\end{bmatrix}
\begin{bmatrix}
\matu_n \\
\matxi_n
\end{bmatrix}
=
\begin{bmatrix}
\mF_{\tau} \\
\mF
\end{bmatrix}
\]
Here, the notations $\matu_n$ and $\matxi_n$ represent the vectors containing the nodal values of the unknown functions. The structure of the matrices $\mK, \mP_{\tau}, \mD$ and $\mM$ can be observed directly from Eqs.~\eqref{disceq1} and \eqref{disceq2}. In this implicit method, additional degrees of freedom are introduced to compute the $L^2$ projection of the FE residual onto the FE space. However, these extra degrees of freedom can be eliminated via static condensation, resulting in
\begin{equation*}
  \left(\mK+\mP_{\tau}\mM^{-1}\mD\right)\matu=\mF_{\tau}+\mP_{\tau}\mM^{-1}\mF
\end{equation*} 
Note that in FE formulations where mass lumping is applied (i.e., with a diagonalizable Gram matrix $\mM$), the computational cost of performing the static condensation is low compared to that of the system of equations. In other applications of OSGS-based stabilization, particularly when solving nonlinear or time-dependent problems, it is more common to compute the projections explicitly, for instance, updating them only after every time step; see, e.g., \cite{Codina2001a,Castanar2020a} and the numerical comparisons in \cite{Gravenkamp2023a,Gravenkamp2023d}.
For the stationary problems considered here, where no initial guess exists, we stick with the implicit approach outlined above.

\section{The dual problem}

\subsection{Strong form of the dual problem}

Consider $Q(u)$ to be a bounded linear functional defined as
\begin{equation*}
  Q:H^1(\Omega)\rightarrow \mathbb{R}
\end{equation*}
yielding a real-valued output in $\mathbb{R}$, which may characterize a physically meaningful quantity of interest. The linear functional $Q(u)$ can also be expressed as 
\begin{equation} \label{linfunc}
  Q(u)= \left\langle q,u\right\rangle
\end{equation}
where $q$ is a representative of $Q$.
The adjoint problem, also known as the dual problem, plays a crucial role in goal-oriented error estimation when using a functional as the quantity of interest. Solving the dual problem helps to identify the regions in the computational domain that contribute most to the functional error. 
The strong form of the dual problem can be expressed as: find $z: \Omega \rightarrow \mathbb{R}$ such that
\begin{align}\label{strdual}
  \begin{split}
      \mathcal{L}^{*}z &= q \qquad \text{in} \quad  \Omega \\
      \hspace{0.3cm}z &= 0 \qquad  \text{on} \quad \Gamma
  \end{split}
\end{align}
where $\mathcal{L}^{*}$ is the adjoint differential operator of $\cal L$ as introduced in \eqref{adj-op}, and $q$ is the forcing function of the dual problem we wish to consider. 

\subsection{Weak form of the dual problem and Galerkin finite element approximation}

In this subsection, we establish the weak formulation of the dual problem. 

The weak form of the adjoint or dual problem of the primal problem \eqref{varform} can be written as: find $z \in V$ such that
\begin{equation} \label{varforadj}
	B^{*}(z,v)=Q(v)\ \ \ \ \ \ \forall v\in V 
\end{equation}
with
\begin{align*}
  B^{*}(z,v) &=\int_{\Omega}{k \nabla z\cdot \nabla v} \ \text{d}\Omega - \int_{\Omega}{a \cdot \nabla z v} \ \text{d}\Omega + \int_{\Omega}{s z v} \ \text{d}\Omega\\
  Q(v) &= \int_{\Omega}{q v} \ \text{d}\Omega
\end{align*}
where $B^{*}$ and $Q$ represent the bilinear form and linear functional of the dual problem, respectively. 

The standard Galerkin approximation of \eqref{varforadj} consist of finding $z_h \in V_h$ such that 
\begin{equation*} 
	B^{*}(z_h,v_h)=Q(v_h)\ \ \ \ \ \ \forall v_h\in V_h 
\end{equation*}

\subsection{Variational multiscale formulation of the dual problem}

Since the Galerkin method fails to provide stable solutions for the convection-dominated problems, the primal problem is discretized using the VMS framework to ensure stability. Similarly, the dual problem is likely to exhibit numerical instabilities using the Galerkin method, and this suggests employing the VMS method for the discretization of the dual problem as well.

We proceed to obtain the stabilized form analogously to the primal problem and repeat the key steps to highlight the crucial differences. First, we decompose the trial functions $z$ and test functions $v$ into coarse and fine scales as
\begin{align*} 
  z&=z_h+z' \qquad z_h \in V_h, ~z' \in V' \\
  v&=v_h+v' \qquad v_h \in V_h, ~v' \in V'
\end{align*}
where $z_h$ and $z'$ represent the coarse and fine-scale components of the dual problem, respectively. Following the above decomposition, the weak form of the dual problem \eqref{varforadj} is split into two subproblems: find $z_h \in V_h$ and $z' \in V'$ such that
\begin{align}
  B^{*}(z_h,v_h)+B^{*}(z',v_h) & = Q(v_h) \qquad \forall v_h\in V_{h} \label{coarsescaleadj} \\ 
  B^{*}(z_h,v')+B^{*}(z',v') & = Q(v') \qquad \forall v'\in V' \label{finescaleadj}
\end{align}
Following the standard procedure of the VMS method, the coarse-scale problem \eqref{coarsescaleadj} can be rewritten using the differential operator to obtain
\begin{equation}
  B^{*}(z_h,v_h)+\sum_{K} {\left\langle z',\mathcal{L} v_h\right\rangle}_K = Q(v_h) \qquad \forall v_h\in V_{h}  \label{dual-stab-fe}
\end{equation}
where we have assumed that $z'$ is negligible on the element boundaries and made use of the fact that $\nabla\cdot a = 0$.
At this point, we require an adequate approximation of the SGS or fine-scale component $z'$. Following the procedure used for the primal problem, the fine-scale problem \eqref{finescaleadj} can be approximated as
\begin{align*}
  B^{*}(z',v') &= Q(v')-B^{*}(z_h,v')\\
           &= Q(v')-\sum_K{{\left\langle {\mathcal{L}^{*}z_h},v'\right\rangle }_K}\nonumber\\
           &= \sum_K{{\left\langle {\mathcal{R}^{*}z_h},v'\right\rangle }_K}\nonumber\\
           & = \sum_K \langle {\cal L}^{*} z',v'\rangle_K\nonumber
\end{align*}
where $\mathcal{R}^{*}z_h = q-\mathcal{L}^{*}z_h$ is the residual of the dual problem. 
Using the same arguments as for the primal problem, the SGS $z'$ can be approximated as 
\begin{equation} \label{sgsapprdual}
  z'  \approx z_h'  := \tau P^{\bot}(\mathcal{R}^{*}z_h) \quad \hbox{(computed element-wise)}
\end{equation}
where the stabilization parameter $\tau$ is the same as for the primal problem, since the difference between $\cal L$ and ${\cal L}^{*}$ is only the sign of the advection velocity of which $\tau$ is independent, Eq.~\eqref{tauelem}.

Using approximation \eqref{sgsapprdual} in \eqref{dual-stab-fe}, the final formulation of the OSGS method for the dual problem can be presented as: find $z_h \in V_h$ such that
\begin{equation} \label{stabformdual}
  B^{*}(z_h,v_h)+\sum_{K} {\left\langle \tau P^{\bot}(\mathcal{R}^{*}z_h),\mathcal{L} v_h\right\rangle}_K = Q(v_h) \qquad \forall v_h\in V_{h} 
\end{equation}
As for the primal problem, the approximation of $z'$ by $z'_h$ in \eqref{sgsapprdual} induces an error, which we now call $\ferrord$. The exact solution of the dual problem is expressed as
\begin{align}
z = z_h + z'_h + \ferrord \label{ferrord}
\end{align}

\subsection{Implementation of the OSGS formulation for the dual problem}

We now present the implementation procedure for the dual problem \eqref{stabformdual}. We have utilized the same implicit methodology as previously applied to the primal problem. We denote again the broken $L^2$ projection of the residual associated with the dual problem onto $V_h$ as follows: 
\begin{equation*}
  \xi_{h,{\rm d}} = P_h(\mathcal{R}^{*}z_h)=P_h(q-\mathcal{L}^{*}z_h)
\end{equation*}
The computation of $\xi_{h,{\rm d}}$ carried out by solving the following equation:
\begin{equation} \label{projcompdual}
  \sum_K \langle \mathcal{L}^{*}z_h, \zeta_h\rangle_K + \sum_{K} {\left\langle \xi_{h,{\rm d}},\zeta_h\right\rangle}_K =\left\langle q,\zeta_h\right\rangle
  \quad \forall \zeta_h\in V_h
\end{equation}
The stabilized formulation of the dual problem can be written as a system of equations: 
\begin{align}
  B^{*}(z_h,v_h) - \sum_K \tau {\left\langle \mathcal{L}^{*}z_h, \mathcal{L}v_h \right\rangle}_K- \sum_K \tau {\left\langle \xi_{h,{\rm d}}, \mathcal{L}v_h \right\rangle}_K &=Q(v_h) - \sum_K \tau {\left\langle q,\mathcal{L}v_h\right\rangle}_K \label{disceq1dual}\\
  \sum_K \langle \mathcal{L}^{*}z_h, \zeta_h\rangle_K + \sum_K {\left\langle \xi_{h,{\rm d}}, \zeta_h\right\rangle}_K &=Q(\zeta_h) \label{disceq2dual}
\end{align}
Analogously to the implementation of the primal problem, the system of equations can be written in the matrix representation as:
\[
\begin{bmatrix}
\mK^\ast & -\mP_{\tau}^\ast \\
\mD^\ast & \mM^\ast
\end{bmatrix}
\begin{bmatrix}
\mz_{n} \\
\xi_{n,{\rm d}}
\end{bmatrix}
=
\begin{bmatrix}
\mQ_{\tau} \\
\mQ
\end{bmatrix}
\]
Here, $\mz_n$ and $\xi_{n,{\rm d}}$ correspond to the arrays of nodal unknowns of the dual problem. Matrices $\mK^\ast, \mP^\ast, \mD^\ast$ and $\mM^\ast$ are defined by comparison with Eqs.~\eqref{disceq1dual} and \eqref{disceq2dual}.

\section{Goal-oriented error estimation}

We continue by introducing a VMS-based goal-oriented a posteriori error estimation framework for the OSGS stabilized method. Our goal is to derive an accurate representation of the error in the quantity of interest, $Q(u-u_h)$. The expressions for the goal-oriented error estimates are derived using the approaches outlined in \cite{granzow2017output,hauke2009variational}, where the VMS method was employed to estimate the error in the quantities of interest by approximating the exact error representation as a function of the SGSs. Our main contribution is twofold. On the one hand, we use the OSGS as the stabilization technique and, on the other hand, we explicitly indicate in our analysis where the approximation of the exact SGSs by the modeled ones is applied. This allows one to track the error incurred in this approximation of the SGSs. 

In the following subsections, we proceed to present two distinct approaches for the OSGS stabilized method based on utilizing the SGS components that emerge from either the primal or the dual problem of the VMS formulations. Also, we prove the equivalence between both approaches. Notably, the error estimators yield identical results globally, although they do not exactly match at the element level. In the subsequent subsections, the notation $\sum_K {\left\langle \cdot,\cdot \right\rangle}_K \equiv\left\langle \cdot,\cdot \right\rangle_h$  will be used to simplify some expressions. 

\subsection{Explicit approach}

In this section, we present a VMS-based explicit approach, where the error is estimated by post-processing the FE solution.  By approximating the SGS component $u'\approx u'_h$, the error in the quantity of interest can be approximated as $Q(u'_h)$, providing a direct and computationally efficient means of error estimation.

Recalling the VMS concept, the decomposition of the continuous solution $u$ as $u=u_h+u'$, along with the linearity of the functional $Q(\cdot)$, leads to the following error representation, assuming we have access to the exact expression of the SGS $u'$:
\begin{equation*}
  Q(u)-Q(u_h)=Q(u')
\end{equation*}
This expression shows that the exact error in the quantity of interest is directly linked with the SGS $u'$, which represents the part not captured by the FE solution and can be considered the exact numerical error. 

However, we do not have at our disposal the exact SGS $u'$ (and consequently the exact projection $u_h$ stemming from the splitting $V = V_h+ V'$).  What we do have is \eqref{ferrorp}, and therefore:
\begin{align*}
  Q(u)-Q(u_h) &= \left\langle q,u\right\rangle-\left\langle q,u_h\right\rangle   &&\text{using} \ \eqref{linfunc} \\  
  &= \left\langle q,u-u_h\right\rangle  &&\text{by linearity}\\  
  &= \left\langle q,u'_h + \ferrorp  \right\rangle   &&\text{using} \ \eqref{ferrorp}
\end{align*}
Our proposal to estimate $\vert Q(u)-Q(u_h) \vert$ is {\em to assume that $u_h' \approx u'$ is a good enough approximation for computing the effect of the SGS on the error in $Q$}, and therefore to neglect the effect of $\ferrorp$. Furthermore, we assume that the integrals involved in the previous expression can be computed element-wise, i.e., replacing $\langle\cdot, \cdot\rangle$ by $\langle\cdot, \cdot\rangle_h$. Thus, the a posteriori error estimate in $Q$ we propose is:
\begin{align}
  Q(u)-Q(u_h)  & \approx  \left\langle q,u'_h  \right\rangle_h  \nonumber \\
 & =   \left\langle q,\tau P^{\bot}(\mathcal{R}u_h)\right\rangle_h   &&\text{using} \ \eqref{subgridOSGS} \label{explicit-st}
 \end{align}
 This expression can also be split into the element contributions to the total error as
 \begin{align}
  Q(u)-Q(u_h)   & \approx \eta_1 = \sum_K \eta_1^K , \label{expli1} \\
 \eta_{1}^{K} & := {\left\langle q,\tau P^{\bot }(\mathcal{R}u_h)\right\rangle}_K \label{expli2} 
\end{align}
$\eta_1$ is the proposed goal-oriented explicit a posteriori error estimator, while $\eta_1^K$ is its element-wise approximation. As $\eta_1$ depends on the solution of the primal problem $u_h$, this approach estimates the error without requiring an additional solution of the dual problem. We have observed numerically that $\eta_1$ performs very well in the problems we have tested. 

\subsection{Implicit approach}

The implicit approach to error estimation relies on the duality principle, a well-established concept in this field. This method requires solving an additional adjoint (dual) problem \eqref{strdual} to obtain error estimates. In contrast to the explicit approach, it estimates the error in the quantity of interest by utilizing the information from both the primal and adjoint solutions. 

The representation of the error in the quantity of interest can be obtained using the exact solution $u$ of the primal problem \eqref{strprim}, the exact solution $z$ of the dual problem \eqref{strdual}, the FE solution $u_h$ of the stabilized formulation of the primal problem \eqref{stabformprim}, and the FE solution $z_h$ of the dual problem \eqref{stabformdual}. It is derived as follows:
\begin{align*}
  Q(u)-Q(u_h) &= \left\langle q,u\right\rangle-\left\langle q,u_h\right\rangle   &&\text{using} \ \eqref{linfunc} \\  
  &= \left\langle \mathcal{L}^{*}z,u\right\rangle-\left\langle \mathcal{L}^{*}z,u_h\right\rangle    &&\text{using} \ \eqref{strdual} \\
  &= \left\langle z,\mathcal{L}u\right\rangle-\left\langle z,\mathcal{L}u_h\right\rangle 
         &&\text{integrating by parts}\\
  &= \left\langle z,f\right\rangle-\left\langle z,\mathcal{L}u_h\right\rangle
        &&\text{using} \ \eqref{strprim} \\
  &= \left\langle z,\mathcal{R}u_h\right\rangle 
      &&\text{residual definition}\\
  &= \left\langle z,\mathcal{R}u_h\right\rangle -\left\langle z_h,\mathcal{R}u_h\right\rangle + \left\langle \mathcal{L}^{*}z_h,\tau P^{\bot}(\mathcal{R}u_h)\right\rangle_h 
     &&\hbox{using \eqref{eqnotstabform}, \eqref{eqnotstabform-bis} with}~P' = P^\bot\\
  &= \left\langle z-z_h,\mathcal{R}u_h\right\rangle + \left\langle \mathcal{L}^{*}z_h,\tau P^{\bot}(\mathcal{R}u_h)\right\rangle_h     &&\text{by linearity}\\
  &= \left\langle z'_h + \ferrord,\mathcal{R}u_h\right\rangle + \left\langle \mathcal{L}^{*}z_h,\tau P^{\bot}(\mathcal{R}u_h)\right\rangle_h    &&\text{using \eqref{ferrord}} 
\end{align*}
Now we introduce similar approximations as for the explicit approach. To estimate $\vert Q(u)-Q(u_h) \vert$, {\em we assume that $z_h' \approx z'$ is a good enough approximation to compute the effect of the SGS on the error in $Q$}, and therefore we neglect the effect of $\ferrord$. Likewise, we also assume that the integrals involved in the previous expression can be computed element-wise. This yields:
\begin{align}
   Q(u)-Q(u_h) & \approx   \left\langle z'_h ,\mathcal{R}u_h\right\rangle_h + \left\langle \mathcal{L}^{*}z_h,\tau P^{\bot}(\mathcal{R}u_h)\right\rangle_h  && \nonumber \\
  &=   \left\langle \tau P^{\bot}(\mathcal{R}^{*}z_h),\mathcal{R}u_h\right\rangle_h+ \left\langle \mathcal{L}^{*}z_h,\tau P^{\bot}(\mathcal{R}u_h)\right\rangle_h     &&\text{using} \ \eqref{sgsapprdual} \label{implicit-st}
\end{align}
This derivation leads to the construction of an error estimator of the form:
\begin{align}
  Q(u)-Q(u_h)  & \approx \eta_2 =  \sum_K \eta_2^K, \label{expli1} \\
 \eta_{2}^{K} & := \left\langle \tau P^{\bot}(\mathcal{R}^{*}z_h),\mathcal{R}u_h\right\rangle_K+ \left\langle \mathcal{L}^{*}z_h,\tau P^{\bot}(\mathcal{R}u_h)\right\rangle_K  \label{expli2} 
\end{align}
The first term accounts for the interaction between the dual SGS component and the primal residual. The second term captures the effect of the primal SGS and the adjoint differential operator. 

\subsection{On the equivalence between error estimators}

As shown in \cite{granzow2017output}, both the explicit and implicit approaches to obtaining the a posteriori error for the quantity of interest are equivalent when the exact expression of the SGSs is used. In our case, considering the approximation of the SGSs we propose, $\eta_1$ and $\eta_2$ will generally be different, but they are expected to be similar. This is due to the fact that:
\begin{align}
 & \left\langle \tau P^{\bot}(\mathcal{R}^{*}z_h),\mathcal{R}u_h\right\rangle_h 
 + \left\langle \mathcal{L}^{*}z_h,\tau P^{\bot}(\mathcal{R}u_h)\right\rangle_h  \nonumber \\
  &\qquad =\left\langle \mathcal{R}^{*}z_h,\tau P^{\bot}(\mathcal{R}u_h)\right\rangle_h+ \left\langle \mathcal{L}^{*}z_h,\tau P^{\bot}(\mathcal{R}u_h)\right\rangle_h       && \nonumber\\
  &\qquad =\left\langle \mathcal{R}^{*}z_h +\mathcal{L}^{*}z_h,\tau P^{\bot}(\mathcal{R}u_h)\right\rangle_h     &&\nonumber\\
  &\qquad =\left\langle q,\tau P^{\bot}(\mathcal{R}u_h)\right\rangle_h       &&\nonumber
\end{align}
where we have used that the broken projection $P^\bot$ satisfies
\begin{align}
\sum_K \langle P^\bot f ,  g \rangle_K = \sum_K \langle P^\bot f ,  P^\bot g \rangle_K = \sum_K \langle f , P^\bot  g \rangle_K\label{prop-proj}
\end{align}
for any functions $f, g \in L^2(K)$ for all $K\in {\cal T}_h$. 

Therefore, the starting points of the explicit error estimate \eqref{explicit-st} and of the implicit one \eqref{implicit-st} are {\em identical}, even if they are obtained by neglecting errors $\ferrorp$ and $\ferrord$, respectively, which are, in principle, different. However, \eqref{prop-proj} holds globally in $\Omega$, but not element-wise. Thus, we should generally expect that $\eta_1^K \not = \eta_2^K$. We will evaluate this difference numerically.

\subsection{A posteriori goal-oriented error in terms of the a posteriori global error}

A classical way to obtain a posteriori goal-oriented error estimates is to rely on a posteriori error estimates in global norms (see, e.g., \cite{diez-pares-huerta-2010}). Let us describe how to apply this idea in our case and explain why it is not particularly effective (unless additional analysis provides deeper insight). 

Let us write the stabilized discrete form of the primal problem~\eqref{stabformprim} and of the dual problem~\eqref{stabformdual} as 
\begin{align}
B_{\rm stab}(u_h, v_h)  & :=
B(u_h,v_h)-\sum_K{{\left\langle \tau P^{\bot}(\mathcal{L}u_h),{\mathcal{L}}^*v_h\right\rangle }_K}\cl
& = L_{\rm stab}(v_h) := L(v_h) -\sum_K{{\left\langle \tau P^{\bot}f,{\mathcal{L}}^*v_h\right\rangle }_K} \quad \forall v_h \in V_h \cl
B^{*}_{\rm stab}(z_h, v_h) & :=
B^{*}(z_h,v_h)-\sum_{K} {\left\langle \tau P^{\bot}(\mathcal{L}^{*}z_h),\mathcal{L} v_h\right\rangle}_K \cl
& = Q_{\rm stab}(v_h) := Q(v_h)  - \sum_K \ {\left\langle \tau P^{\bot}q ,\mathcal{L} v_h\right\rangle}_K  \quad \forall v_h\in V_{h}  \el
\end{align}
Using \eqref{prop-proj}, it is clear that $B^{*}_{\rm stab}(z_h,v_h) = B_{\rm stab}(v_h,z_h)$. Furthermore, both the primal and the dual problems are consistent. In particular, if $z$ solves~\eqref{strdual}, then $ B^{*}_{\rm stab}(z, v) = Q_{\rm stab}(v) $ for all $v\in V$ locally regular enough. In view of these observations, we have that
\begin{alignat}{3}
Q_{\rm stab}(u) - Q_{\rm stab}(u_h) & = Q_{\rm stab}(u-u_h) && \quad \hbox{$Q_{\rm stab}$ is linear} \cl
& = B_{\rm stab}^{*}(z,u-u_h) && \quad \hbox{the stabilized dual problem is consistent} \cl
& = B_{\rm stab}(u-u_h , z )  && \quad \hbox{}  \label{eq:rm1} \\
& = B_{\rm stab}(u-u_h , z - z_h )  && \quad \hbox{the stabilized primal problem is consistent}  \label{qoi-global}
\end{alignat}
In cases in which $B_{\rm stab}$ is symmetric and positive (semi) definite, it induces a (semi) norm $\Vert \cdot \Vert_{B_{\rm stab}}$, i.e., $B_{\rm stab}(v ,v ) = \Vert v\Vert^2_{B_{\rm stab}}$. The previous inequality allows one to obtain upper and lower estimates for $\vert Q_{\rm stab}(u) - Q_{\rm stab}(u_h)\vert $ in terms of estimates for $\Vert u - u_h \Vert_{B_{\rm stab}}$ and $\Vert z - z_h \Vert_{B_{\rm stab}}$. However, in the case of the CDR equation, $B$ is not symmetric, and it does not induce a norm. If one can show that $B_{\rm stab}(u,v) \leq C \Vert u \Vert_{B_1}\Vert v \Vert_{B_2}$ for a certain constant $C>0$ and norms $\Vert \cdot \Vert_{B_1}$ and $\Vert \cdot \Vert_{B_2}$, then
\begin{align*}
\vert Q_{\rm stab}(u) - Q_{\rm stab}(u_h)\vert \leq C \Vert u - u_h \Vert_{B_1}\Vert z - z_h \Vert_{B_2}
\end{align*}
and thus estimating the error for the quantity of interest reduces to estimating a posteriori the error of the primal problem in the global norm $\Vert \cdot \Vert_{B_1}$ and of the dual problem in the global norm $\Vert \cdot \Vert_{B_2}$. However, for the problem we consider, as it has already been said such norms are either not robust with respect to convection or, to our knowledge, there are no available a posteriori error estimates for them. 

For CDR problems, the analysis showing that suitable norms exist {\em using the Galerkin method} can be found in~\cite{verfurth2005robust}. Let $\Vert \cdot \Vert$ be the $L^2$ norm. If we define
\begin{align*}
\Vert v \Vert_{B_2}^2 := k \Vert \nabla v \Vert^2 + s \Vert v \Vert^2, \quad  
\Vert v \Vert_{B_1} := \Vert v \Vert_{B_2} + \sup_{w\in V} \frac{\langle a \cdot \nabla v , w\rangle}{\Vert v \Vert_{B_2}} 
\end{align*}
then $ \vert Q(u) - Q(u_h)\vert \leq C \Vert u - u_h \Vert_{B_1}\Vert z - z_h \Vert_{B_2}$, and there are a posteriori estimates for  $\Vert u - u_h \Vert_{B_1}$ and $\Vert z - z_h \Vert_{B_2}$, see~\cite{verfurth2005robust}. However, Galerkin solutions are unstable when convection dominates, and we need $u_h$ to be a solution of the Galerkin problem when deriving~\eqref{qoi-global} (although not to obtain an error estimate for  $\Vert u - u_h \Vert_{B_1}$, see~\cite{tobiska2015robust}). Furthermore, the norm $\Vert \cdot \Vert_{B_1}$ is rather weak. We have obtained global a posteriori error estimates in stronger norms (see \cite{codina-gravenkamp-khan-2025}) but under unprovable assumptions. 

Let us also point out that Eq.~\eqref{eq:rm1} is the basis of the dual weighted residual (DWR) approach (see the review in \cite{becker-rannacher-2001}), which consists in replacing $z$ by an approximation $z_{h'}$, usually obtained with a finer approximation than $z_h$. This approach is used for example in \cite{duprez-et-al-2020} in the context of goal-oriented error estimation in solid mechanics applied to soft tissues. The difficulty of applying this strategy in our case is again the need of an adequate norm of $u-u_h$.

This discussion justifies the approach followed in the previous subsections to obtain a posteriori error estimates for quantities of interest without relying on a posteriori global error estimates.

\section{Numerical Results}

This section presents numerical experiments that evaluate the performance of the proposed goal-oriented error estimates for CDR problems in one and two dimensions, all cases being dominated by convection, that is the situation in which we are interested. For consistency, we employed the same stabilization scheme, namely, the OSGS method, for both the primal \eqref{strprim} and dual \eqref{strdual} problems. The stabilized formulations of both are given by Eq.~\eqref{stabformprim} and Eq.~\eqref{stabformdual}, respectively. We have not tested here the more classical ASGS approach (for a comparison between OSGS and ASGS in global a posteriori error estimation, see~\cite{codina-gravenkamp-khan-2025}).

The numerical discretization is performed using the same basis functions and meshes. Moreover, we employed linear interpolations for both the primal and dual problems in all test cases. A goal-oriented error estimation is performed for a specific quantity of interest defined by a linear functional $Q(u)$ \eqref{linfunc}. In all two-dimensional examples, we used uniform bilinear quadrilateral elements for the spatial discretization. These numerical tests demonstrate:
\begin{enumerate}
  \item The performance of the error estimators for a 1D convection-dominated problem.
  \item The accurate recovery of the error on a 2D mesh.
  \item The assessment of the error estimators for a strong boundary layer problem with large gradients.
  \item The error performance in an L-shaped benchmark problem.
\end{enumerate}
The linear functional quantity $Q$ associated with the FE approximation $u_h$ is given by
\begin{align*}
  Q(u_h) &= \left\langle q,u_h\right\rangle
  =\int_{\Omega}qu_h \ \text{d}\Omega \nonumber
\end{align*}
Or equivalently, it can be evaluated element-wise as 
\begin{equation*}
  Q(u_h) = \sum_{K} \int_{K}qu_h \ \text{d}\Omega
\end{equation*}
The performance of the error estimators is evaluated using the effectivity index, defined as
\begin{equation}\label{effind}
  \mathcal{I}_{\mathrm{eff}} = \frac{|Q(u)-Q(u_h)|}{|\eta|}
\end{equation}
Here, the effectivity index $\mathcal{I}_{\mathrm{eff}}$ represents the ratio of the exact error to the estimated error $\eta$ (either $\eta_1$ or $\eta_2$), which should ideally be near or equal to $1$. The global error estimator is computed as 
\begin{equation*}
  \eta_i = \sum_{K} \eta_i^{K},\quad i = 1,2
\end{equation*}
where $\eta_i^{K}$ denotes the local error contribution from element $K$. Furthermore, we describe the behavior of global and local error estimators $\eta_1$ and $\eta_2$. Both error estimators show identical behavior on a global scale, i.e., $\eta = \eta_1 = \eta_2$; however, on the element level, they do not match exactly. In the numerical tests, we analyze the performance of the error estimates in terms of the global effectivity index. In addition, we compare the local error contributions of the explicit and implicit error estimators, $\eta_1^K$ and $\eta_2^K$ for each element $K$, respectively.

\subsection{One-dimensional example} \label{example-1}

\begin{figure}[h]
  \centering
  \subfloat[]{\includegraphics[width =0.43\textwidth]{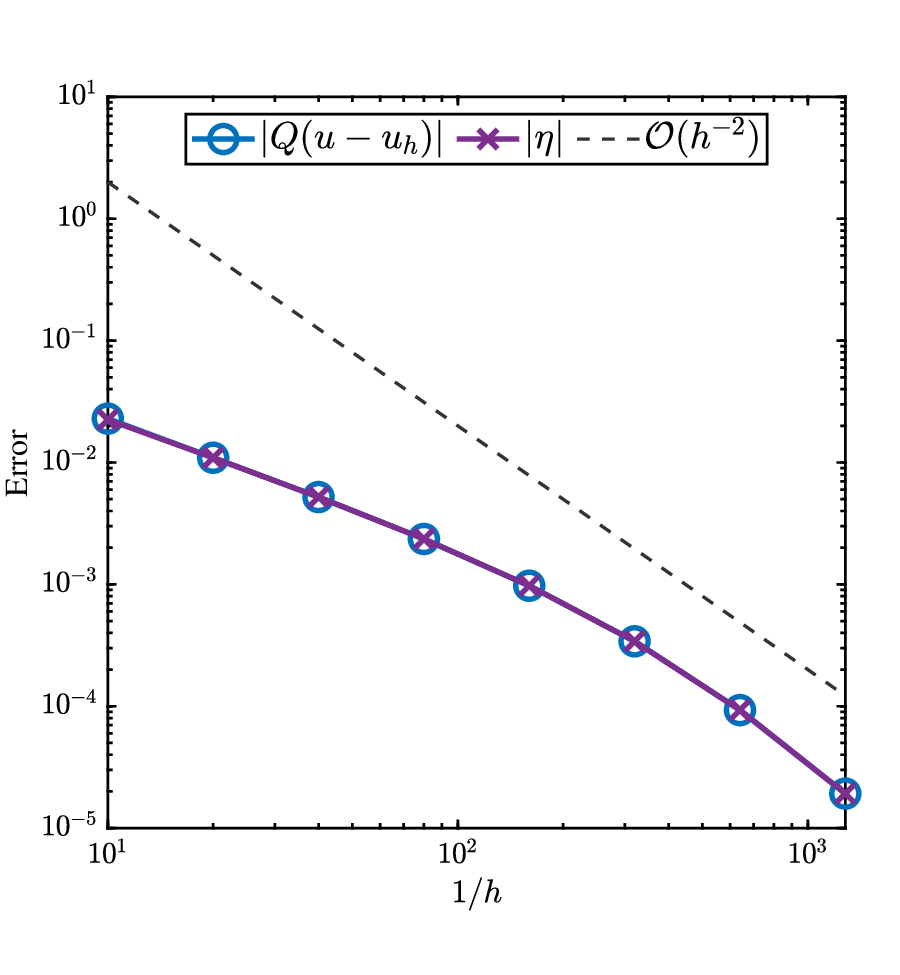}}
  \subfloat[]{\includegraphics[width =0.43\textwidth]{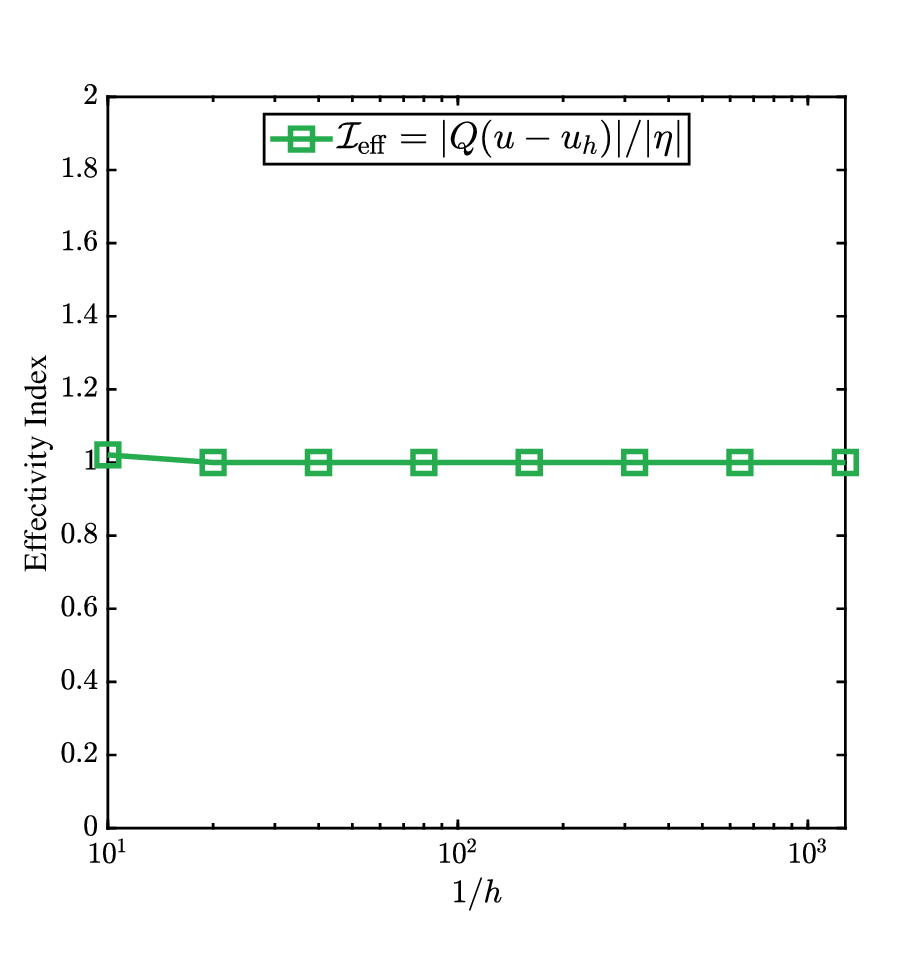}}
  \hfill
  \caption{Results of example \ref{example-1}: Error convergence in the quantity of interest functional $Q$ compared to the goal-oriented error estimator (a);  global effectivity index (b).}
  \label{ex1:errcon-eff}
\end{figure}

The first example is designed to evaluate the performance of the goal-oriented error estimation in terms of the quantity of interest for a 1D convection-dominated problem. For simplicity, we choose the diffusion coefficient as $k=1$ and the reaction coefficient as $s=0.1$ with large convection $a=1000$; hence, the stabilization mechanism is essential. We consider the computational domain
\begin{equation*}
  \Omega =\left\{ x \in \mathbb{R}\ \,|\ \, 0 < x < 1 \right\}
\end{equation*}
where we impose Dirichlet boundary conditions $u(0)=1$ and $u(1)=0$ and assume a vanishing boady load $f=0$. The solution is given by
\begin{equation*}
  u (x)=\frac{1-e^{-1000(1-x)}}{1-e^{-1000}}
\end{equation*}
In this example, we set $q=1$ within the entire domain $\Omega$ in the linear functional \eqref{linfunc}, resulting in an exact value of $Q(u)=0.9990$. We assess the performance of the proposed error estimates in terms of the effectivity index defined in \eqref{effind}.  Fig.~\ref{ex1:errcon-eff} shows the convergence of the error in the global quantity of interest functional $Q$ and the global error estimator $\eta$ for different mesh sizes, where $h$ denotes the element size. In Fig.~\ref{ex1:errcon-eff} (a), a reference line with a convergence rate of $h^{-2}$ is plotted for comparison with the computed error convergence. It can be observed that the error in the quantity of interest functional converges at the same rate as the error estimator. As the global estimators $\eta_1$ and $\eta_2$ produce identical results, i.e, $\eta = \eta_1 = \eta_2$, we plot only one of them to avoid repetition. Fig.~\ref{ex1:errcon-eff} (b) shows that the effectivity index converges to precisely $1.0$ for every mesh. This confirms the accurate recovery of the error for the one-dimensional convection-dominated problem. It also demonstrates that the global functional quantity $Q(u_h)$ converges to the true value of $Q(u)$ as $h \rightarrow 0$. 

\subsection{Two-dimensional example} \label{example-2}
\begin{figure}[h]
  \centering
  \subfloat[]{\includegraphics[width =0.43\textwidth]{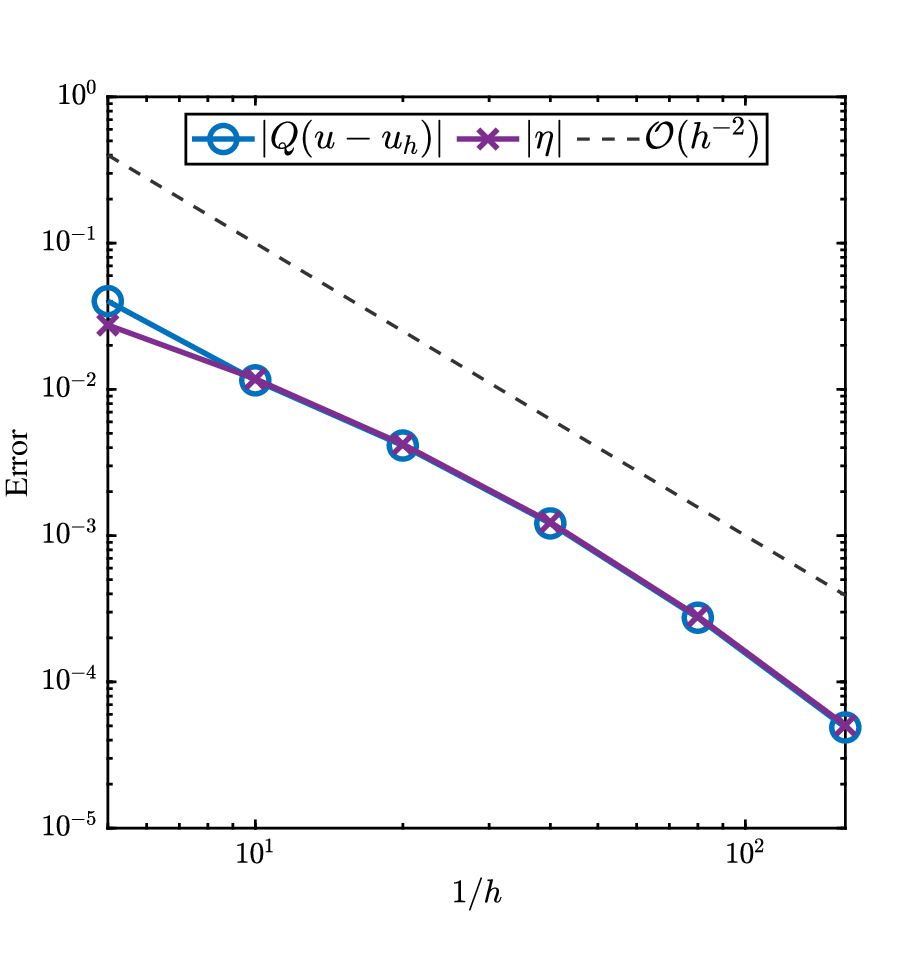}}
  \subfloat[]{\includegraphics[width=0.43\textwidth]{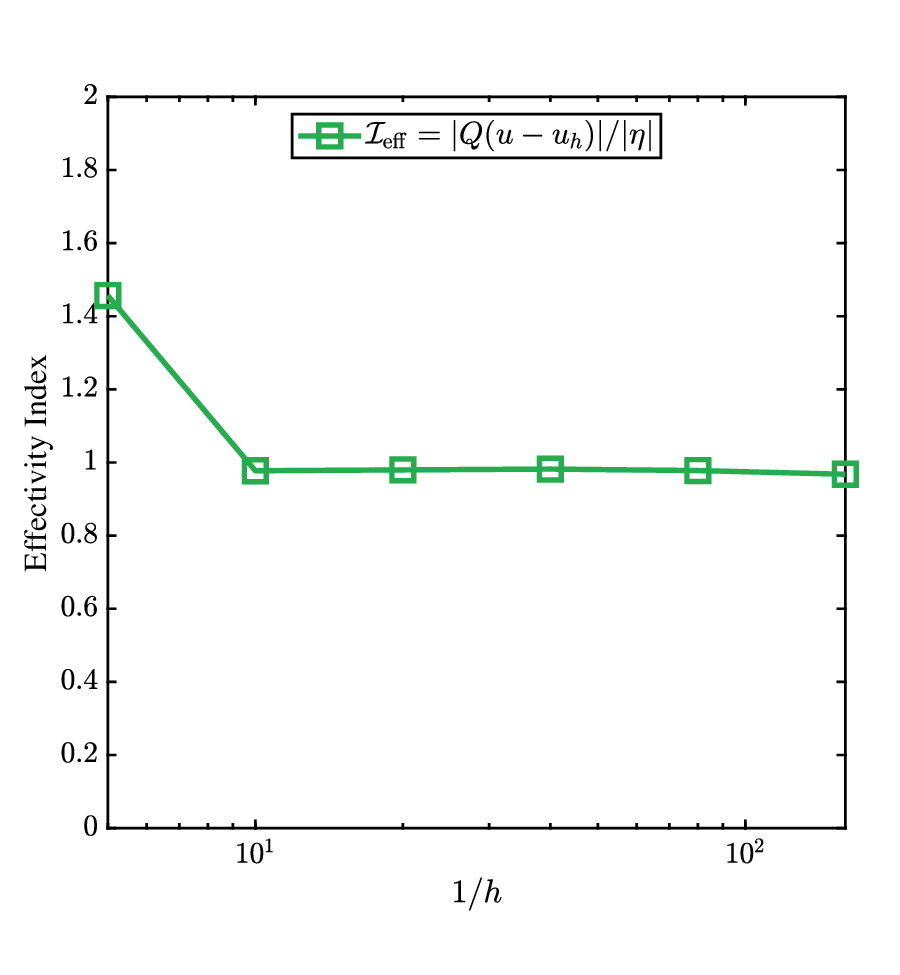}}
  \hfill
  \caption{Results of example \ref{example-2}: Error convergence in the quantity of interest $Q$ compared to the goal-oriented error estimator (a);  global effectivity index (b).}  
  \label{ex2:errconeff}
\end{figure}

\begin{figure}[h]
  \centering
  \subfloat[]{\includegraphics[width =0.43\textwidth]{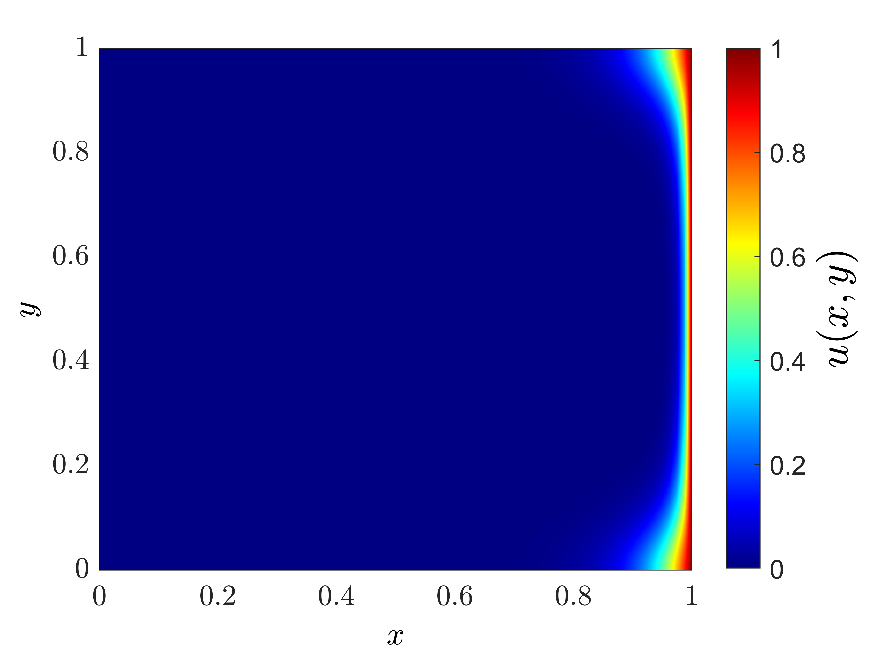}}
  \subfloat[]{\includegraphics[width =0.43\textwidth]{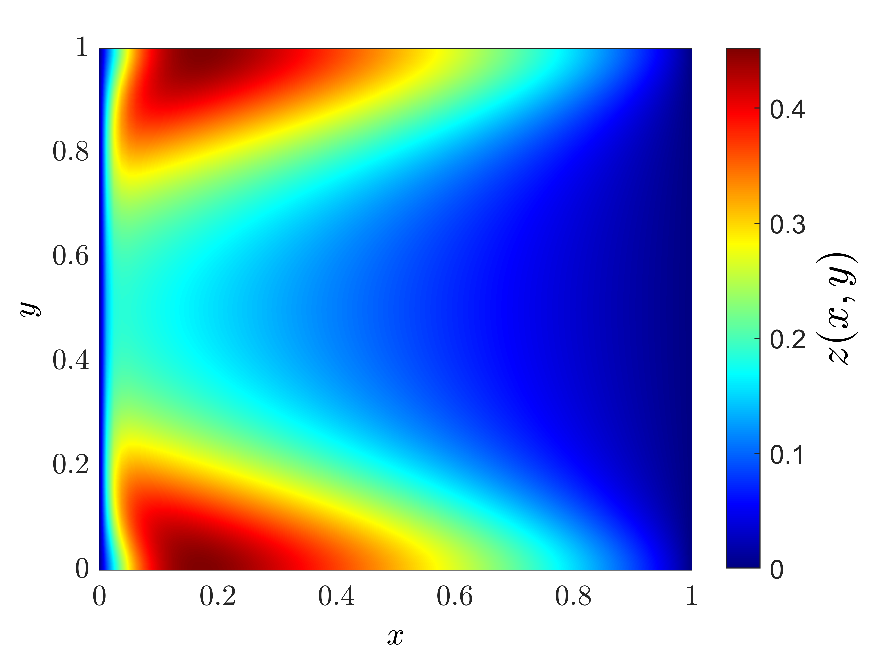}}
  \hfill
  \caption{Results of example \ref{example-2}: (a): Numerical solution of the primal problem \eqref{strprim}. (b): Numerical solution of the dual problem \eqref{strdual}.}
  \label{ex2sol}
\end{figure}

\begin{figure}[h]
  \centering
  \subfloat[]{\includegraphics[width=0.43\textwidth]{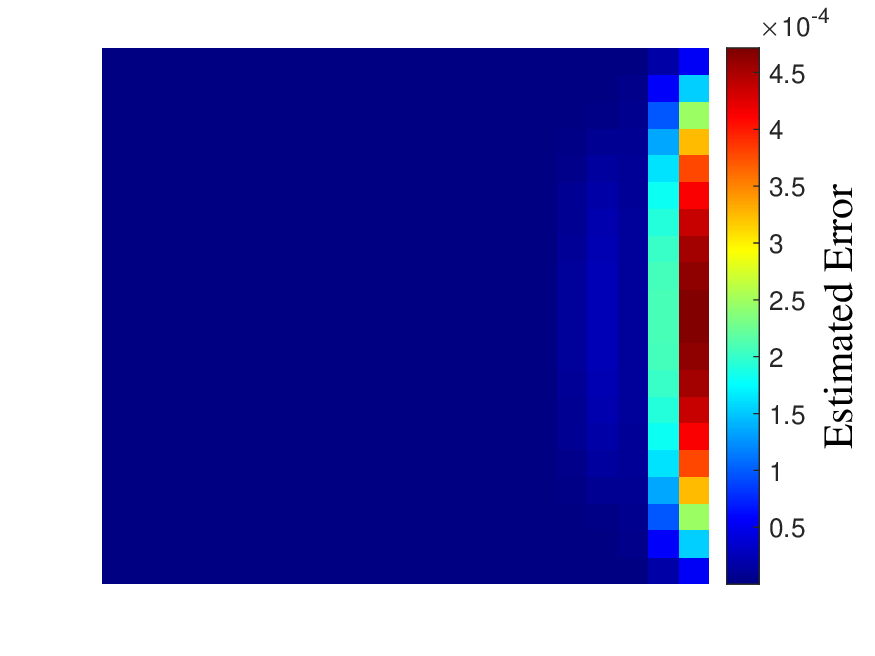}}
  \subfloat[]{\includegraphics[width =0.43\textwidth]{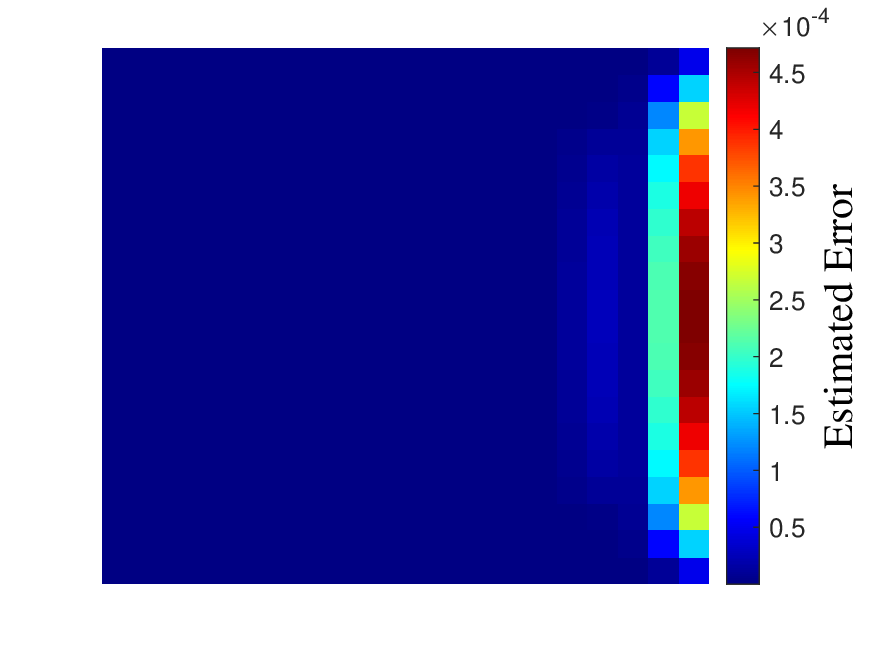}}
  \hfill
  \caption{Results of example \ref{example-2}: Local contribution of the explicit error estimator $\eta_1$ (a) and implicit error estimator $\eta_2$ (b).}
  \label{ex2error}
\end{figure}

In our second example, we test the performance of the goal-oriented error estimators for a two-dimensional problem in a square domain, given by
\begin{equation*}
  \Omega =\left\{ (x,y) \in \mathbb{R}^2\hspace{0.2cm}|\hspace{0.2cm}0 < x 
  < 1, \ 0 < y < 1 \right\} = (0,1)^2
\end{equation*}
We choose the diffusion coefficient $k=0.05$, the reaction coefficient $s=0$, and the convection field  
\begin{equation*}
  a(x,y)= \big(20y(1-y),\,0\big)
\end{equation*}
representing a horizontal parabolic flow profile with a maximum magnitude at $y=0.5$, inducing a horizontal transport-dominated flow from left to right. Dirichlet conditions are prescribed as $u(0,y)=0$ and $u(1,y)=1$ at the inflow and outflow boundaries, respectively. This setup creates a boundary layer near the outflow boundary at $x=1$. We assume again a vanishing forcing term, i.e., $f=0$. Since we do not know the exact solution of the problem, a uniform fine mesh of $640 \times 640$ elements is employed to compute the reference solution, $u_{\rm{ref}}$, which is assumed to be a sufficiently accurate approximation of the exact solution. In this example, the quantity of interest functional is defined as:
\begin{equation*}
  Q(u)= \int_{\Omega} \cos \left( \frac{\uppi x}{5} \right) u \ \text{d}\Omega
\end{equation*}
 A numerical approximation of the exact solution for the chosen quantity of interest,  computed on a sufficiently fine mesh, is given by $Q(u_{\rm{ref}}) \approx Q(u) = 0.0175$.
 Fig.~\ref{ex2:errconeff} shows the computed error in the quantity of interest $Q(u-u_h)$ and the effectivity index for different mesh sizes. It can be observed that the effectivity index converges to a value very close to $1$, demonstrating that the error estimator has accurately captured the error. It also shows that the error in the quantity of interest and the error estimator converge at the same rate as $h \rightarrow 0$. An illustration of the numerical solutions of the primal and dual problems is presented in Fig.~\ref{ex2sol}. The solution is obtained using the OSGS method on an $80 \times 80$ mesh consisting of bilinear quadrilaterals. Fig.~\ref{ex2error} illustrates the local contributions of both the explicit and implicit error estimators. For the results in this figure, we have utilized a uniform mesh of $20 \times 20$ bilinear quadrilaterals. We observe that the local error contribution of the explicit and implicit error estimators is very similar.  In addition, higher error concentrations in the boundary layer region near the outflow boundary at $x=1$ can be observed.

\subsection{A strong boundary layer problem} \label{example-3}

This numerical example focuses on analyzing the behavior of the error estimators for a problem involving a strong boundary layer with large gradients on the unit square, i.e., the same domain as in the previous example. We used $k=0.01$, $a=[1,1]$, $s=10^{-4}$, and the forcing term $f$ on the right-hand side of \eqref{strprim} is chosen such that 
\begin{equation*}
  u(x,y)=\left(x-\frac{1-e^{100x}}{1-e^{100}}\right)\left(y-\frac{1-e^{100y}}{1-e^{100}}\right)
\end{equation*}
is the exact solution. We have imposed homogeneous Dirichlet boundary conditions.  

\begin{figure}[h]
  \centering
  \subfloat[]{\includegraphics[width =0.43\textwidth]{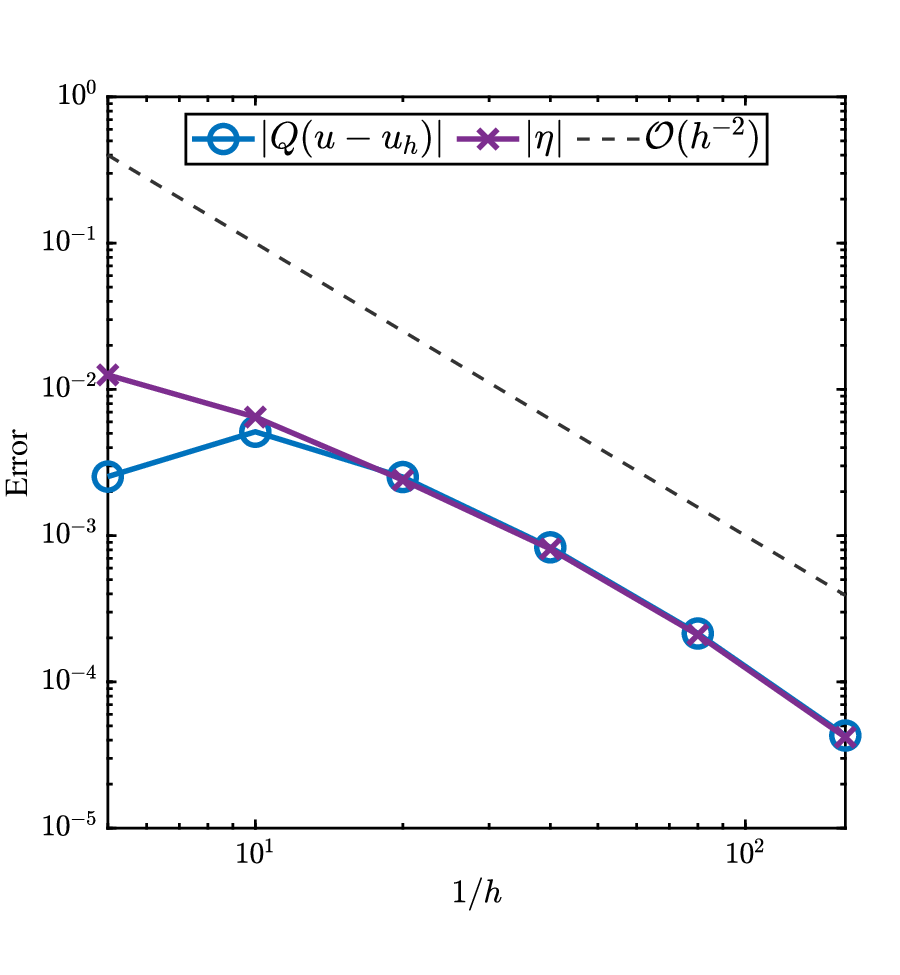}}
  \subfloat[]{\includegraphics[width=0.43\textwidth]{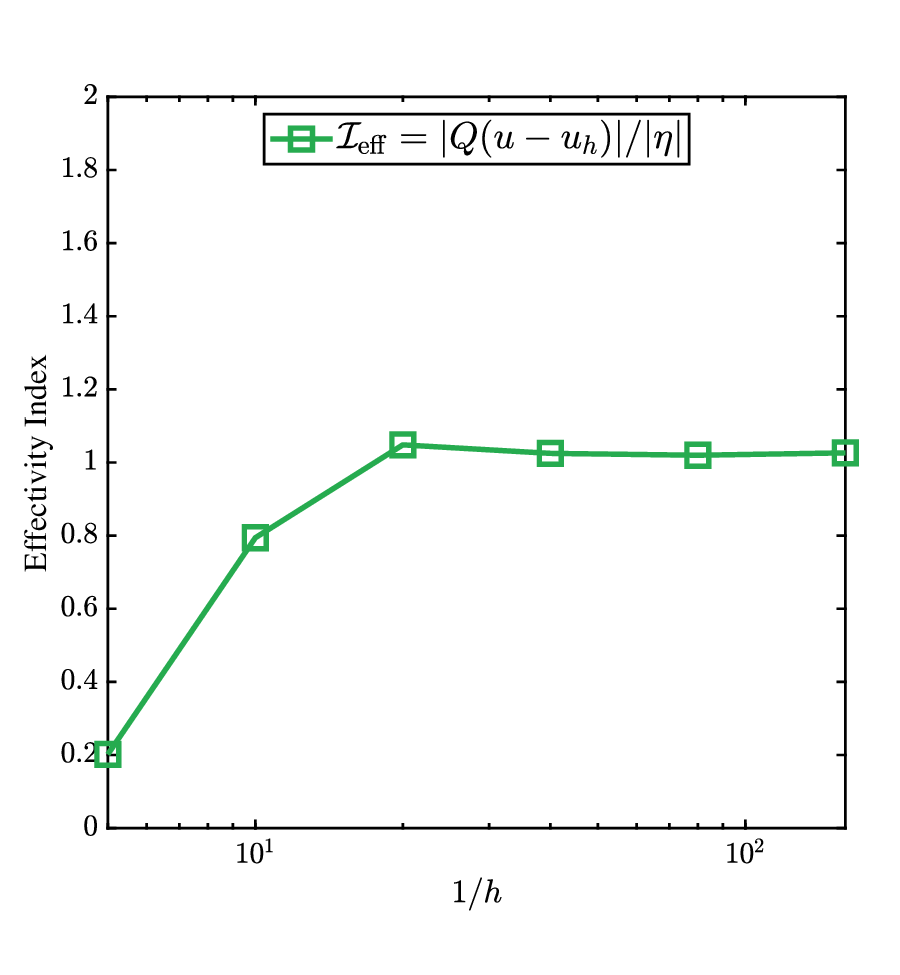}}
  \hfill
  \caption{Results of example \ref{example-3}: Error convergence in the quantity of interest $Q$ compared to the goal-oriented error estimator for the strong boundary layer problem (a);  global effectivity index (b).}
  \label{ex3:errconeff}
\end{figure}

In this example, we are interested in measuring the error in a specific region of interest in a subdomain $\Omega_s=(0.75,1)^2$. In that case, the exact solution for the chosen quantity of interest functional \eqref{linfunc} is $Q(u)=0.0436$ in $\Omega_s$, where $q$ is defined as
\begin{align} 
  \begin{cases}
      q &= 1 \ \ \text{if} \ 0.75 \leq x \leq 1 \ \text{and} \ 0.75 \leq y \leq 1\\
      q &= 0 \ \ \text{otherwise}
  \end{cases}
\end{align}
Fig.~\ref{ex3:errconeff} shows the error convergence plots and the effectivity index for different mesh sizes. Initially, on coarse meshes, the effectivity index deviates from 1. As the mesh is uniformly refined, the error estimates converge to the exact error, resulting in effectivity indices that are very close to unity. This confirms that the error estimators have correctly captured the error and, thus, the linear functional $Q(u_h)$ associated with $u_h$ precisely recovers the exact quantity $Q(u)$. 

\begin{figure}[h]
  \centering
  \subfloat[]{\includegraphics[width=0.43\textwidth]{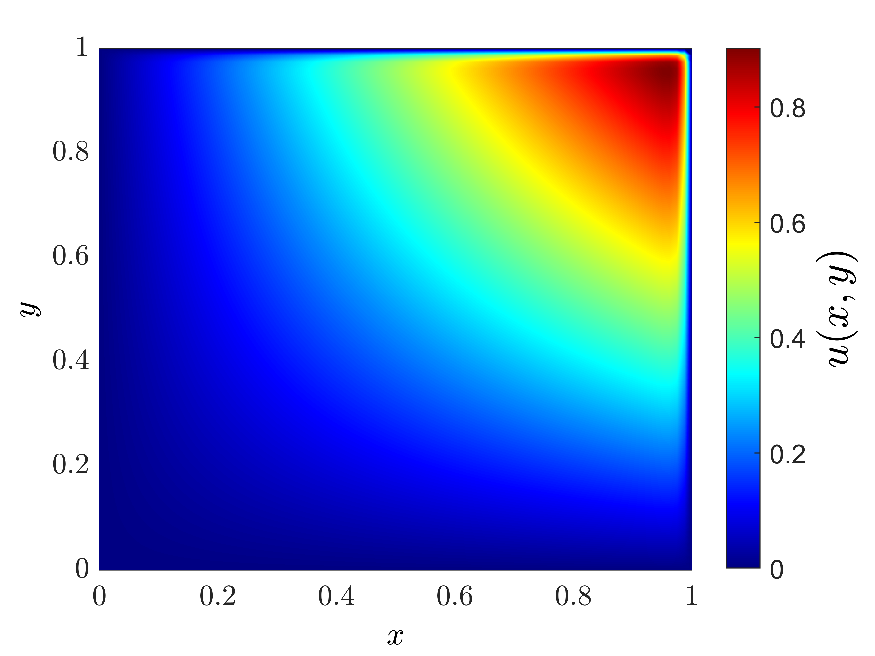}}
  \subfloat[]{\includegraphics[width=0.43\textwidth]{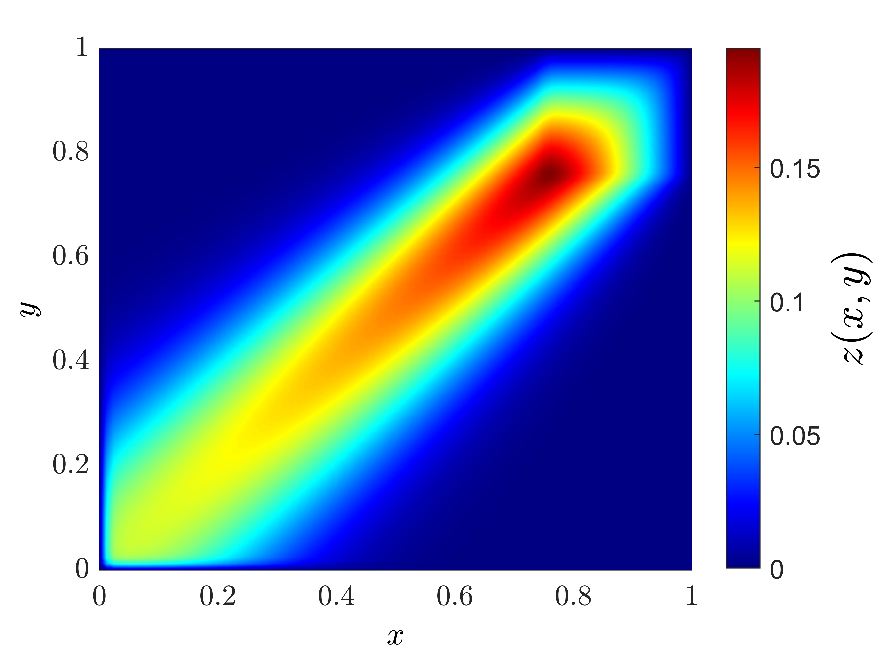}}
  
  \hfill
  \caption{Results of example \ref{example-3}: (a): Finite element solution of the primal problem \eqref{strprim}. (b): Finite element solution of the dual problem \eqref{strdual}.}
  \label{ex3sol}
\end{figure}

Fig.~\ref{ex3sol} (a) shows the numerical solution of the primal problem. This solution exhibits strong boundary layers along $x=1$ and $y=1$. As shown in Fig.~\ref{ex3sol} (b), the dual solution exhibits boundary layers similar to the primal solution. However, due to the reversed convection direction, the layers shifted to the opposite side of the domain compared with the primal problem. These results have been obtained using a uniform mesh of $80 \times 80$ quadrilateral elements.

\begin{figure}[h]
  \centering
  \subfloat[]{\includegraphics[width=0.43 \textwidth]{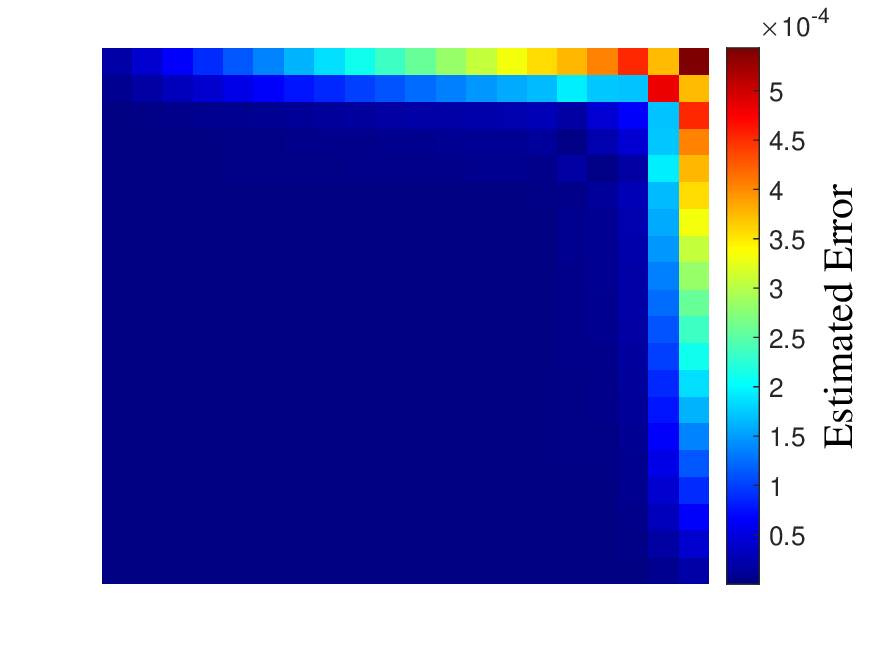}}
  \subfloat[]{\includegraphics[width=0.43 \textwidth]{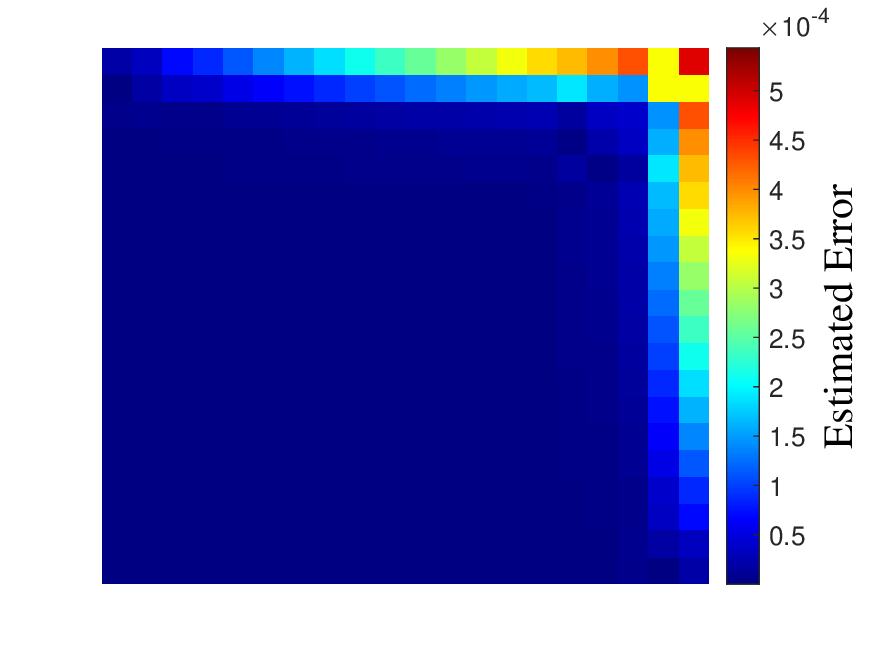}}
  \hfill
  \caption{Results of example \ref{example-3}: Local contribution of the explicit error estimator $\eta_1$ (a) and implicit error estimator $\eta_2$ (b).}
  \label{ex3eta}
\end{figure}

It is useful to study the behavior of the local error estimator in each element. To observe the influence of the error over the entire square domain $\Omega=(0,1)^2$, we choose $q=1$ in $\Omega$. Following this, Fig.~\ref{ex3eta} presents a comparison of the local error estimator contributions to $\eta_1$ and $\eta_2$ (i.e., $\eta^K_1$ and $\eta^K_2$ for all elements $K$). It also reveals that the error is most prominent in the boundary layer region, where the solution changes abruptly.


\subsection{L-shaped domain} \label{example-4}

In this test, we present a convection-dominated problem in an L-shaped domain. We choose
\begin{equation*}
  \Omega = \left(0,1\right) \times \left(0,1\right)\backslash\left[0.5,1\right]\times\left[0,0.5\right].
\end{equation*}
The domain is again discretized using uniform bilinear quadrilateral elements. We chose the diffusion coefficient $k=0.001$, the convection velocity $a=[1,1]$, the reaction coefficient $s=0$, and the forcing function $f=1$. Since the analytical solution of the problem is unknown, a fine mesh consisting of $196608~(= 0.75 \cdot 512^2)$ elements is used to compute the reference solution $u_{\rm ref}$, which is assumed to be a close approximation of the exact solution. For $q=1$ in $\Omega$, the numerical approximation of the exact solution of the chosen quantity of interest, computed on a significantly fine mesh, turns out to be $Q(u_{\rm ref}) \approx Q(u) = 0.2063$. 

\begin{figure}[h]
  \centering
  \subfloat[]{\includegraphics[width =0.43\textwidth]{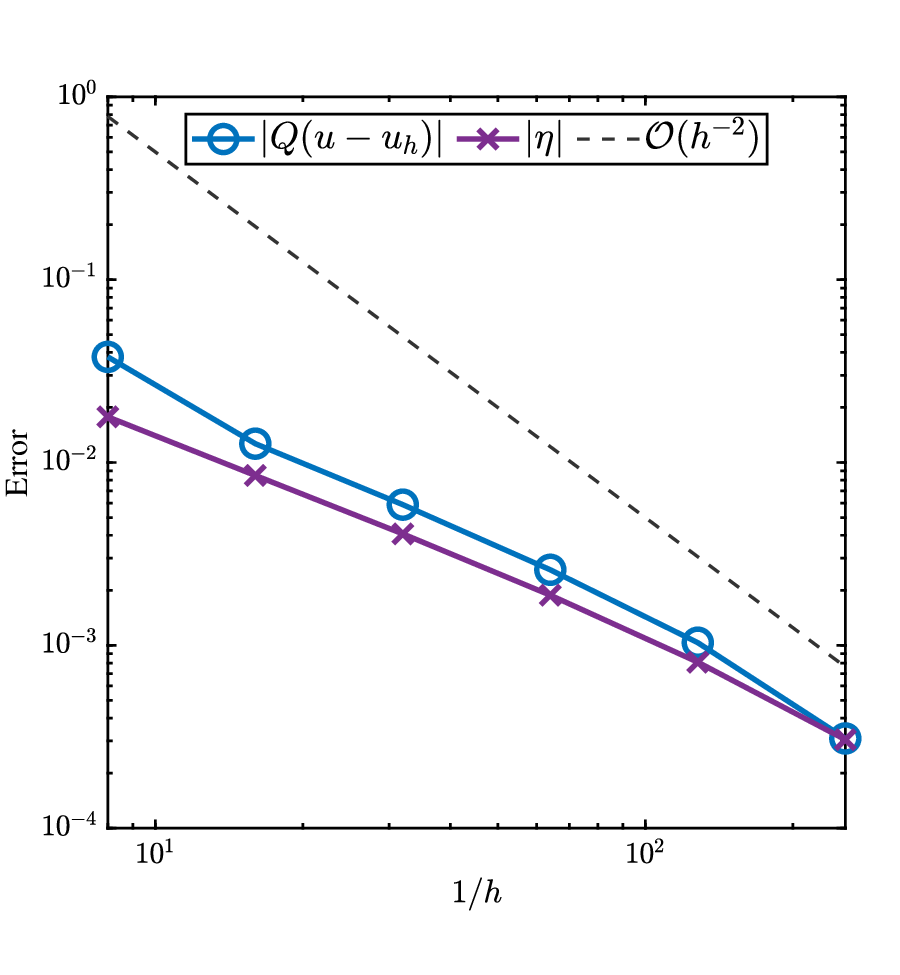}}
  \subfloat[]{\includegraphics[width=0.43\textwidth]{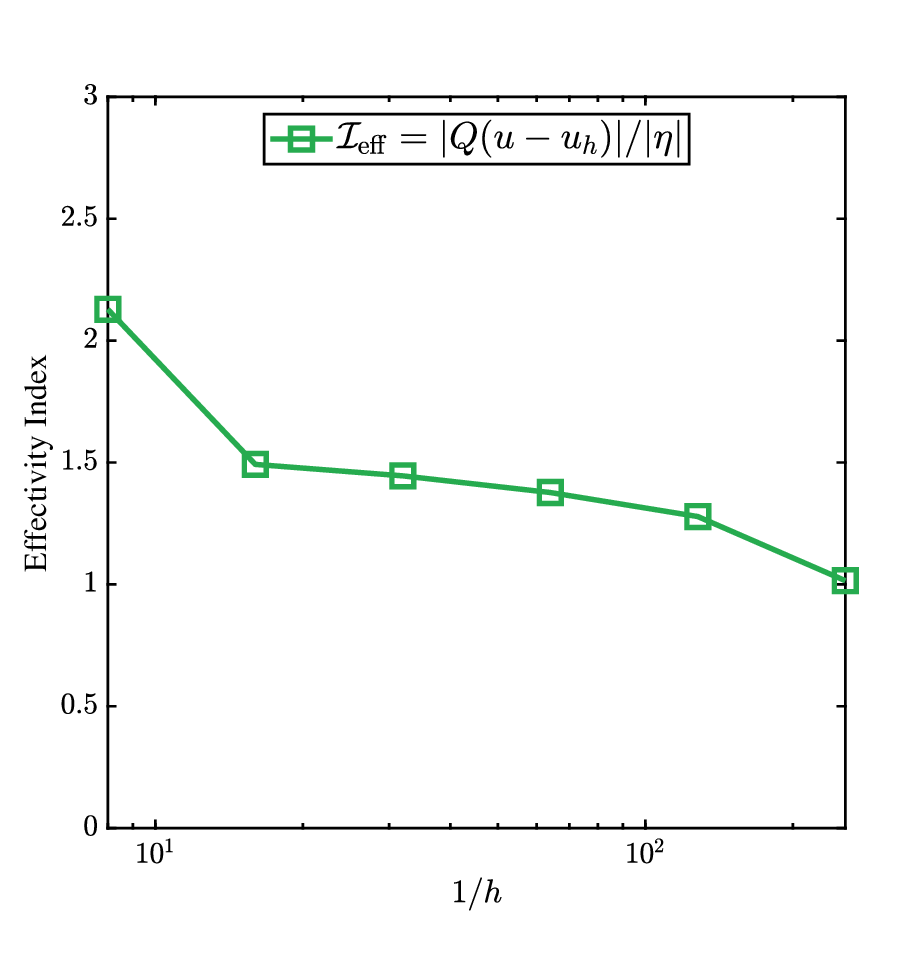}}
  \hfill
  \caption{Results of example \ref{example-4}: Error convergence in the quantity of interest $Q$ compared to the goal-oriented error estimator (a);  global effectivity index (b).}
  \label{ex4:errconeff}
\end{figure}

\begin{figure}[h]
  \centering
  \subfloat[]{\includegraphics[width=0.43 \textwidth]{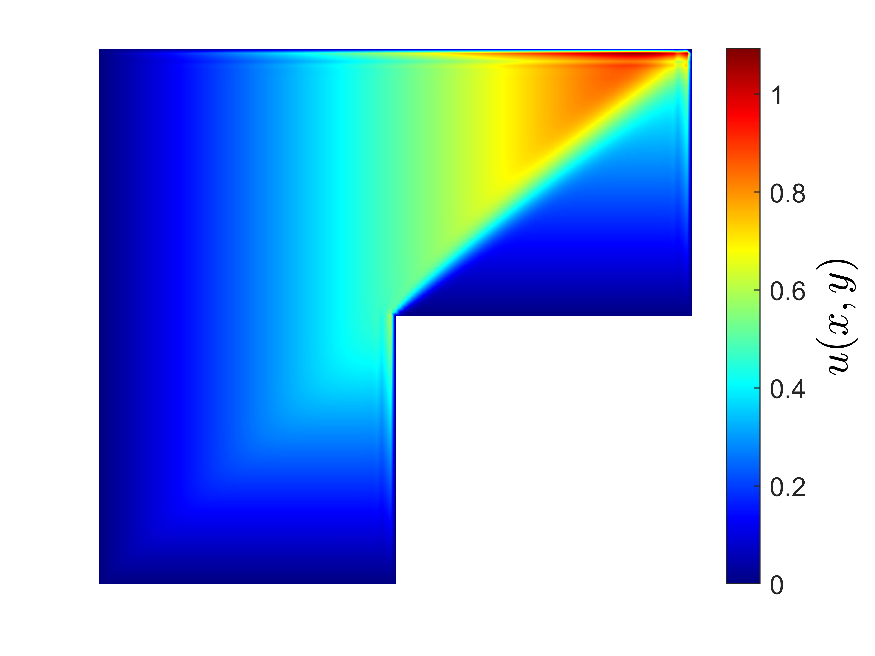}}
  \subfloat[]{\includegraphics[width=0.43 \textwidth]{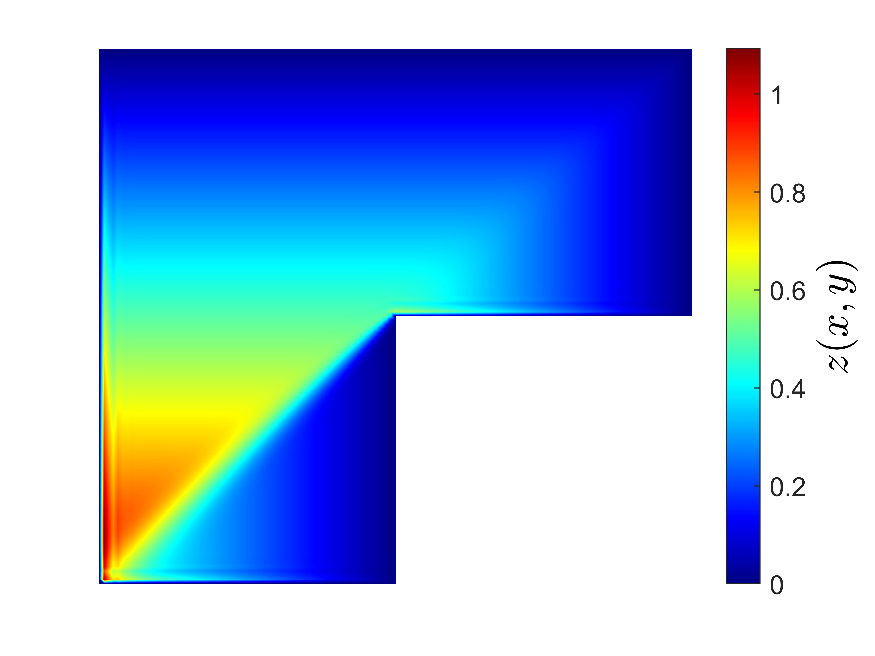}}
  \hfill
  \caption{Results of example \ref{example-4}: Numerical solutions of the primal problem (a) and the dual problem (b).}
  \label{ex4sol}
\end{figure}

\begin{figure}[h]
  \centering
  \subfloat[]{\includegraphics[width=0.43 \textwidth]{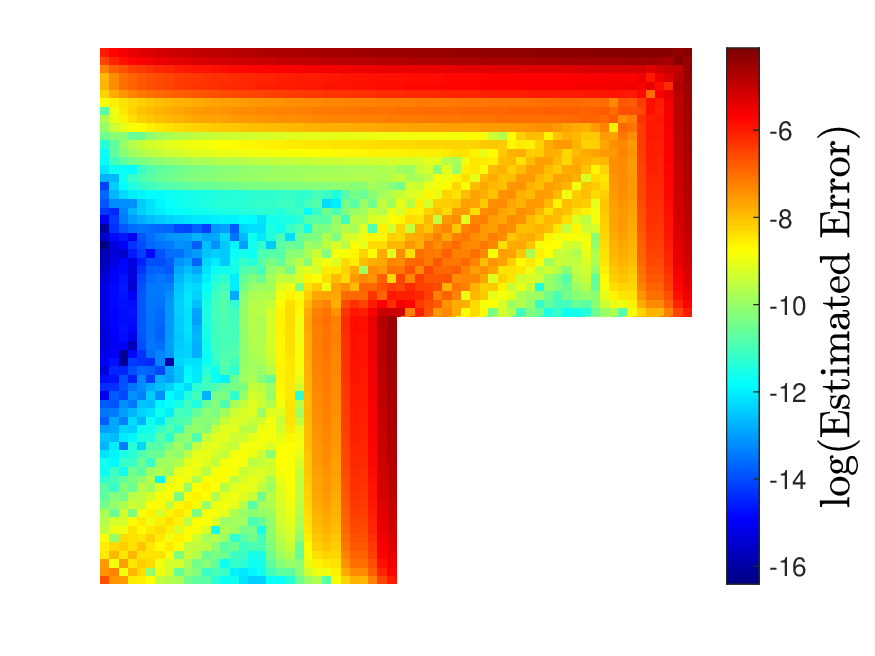}}
  \subfloat[]{\includegraphics[width=0.43 \textwidth]{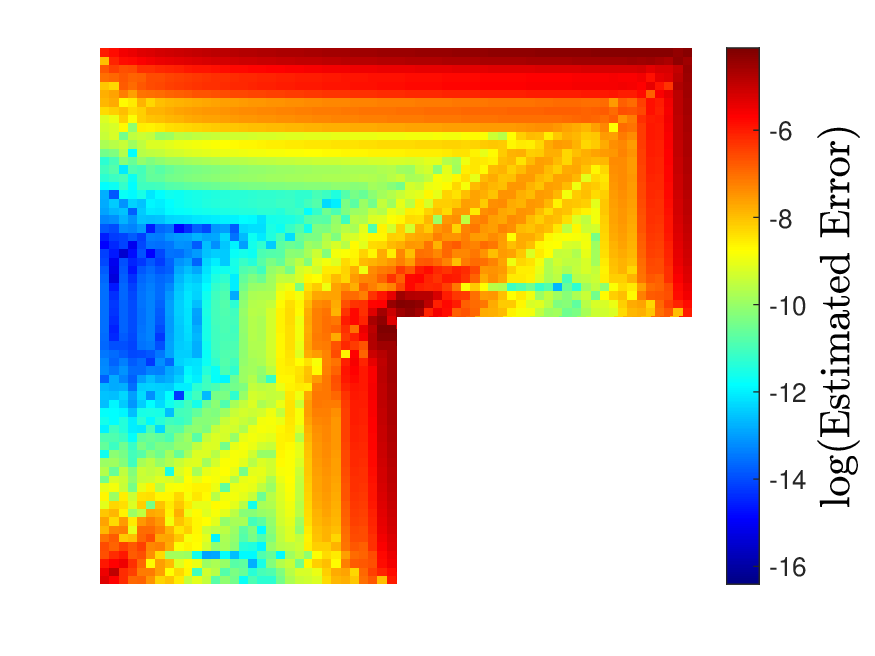}}
  \hfill
  \caption{Results of example \ref{example-4}: Local contribution of the explicit error estimator $\eta_1$ (a) and implicit error estimator $\eta_2$ (b).}
  \label{ex4eta}
\end{figure}

Fig.~\ref{ex4:errconeff} shows the convergence of the error, error estimator, and effectivity indices using uniform mesh refinements. The effectivity indices converge close to $1$ as $h \rightarrow 0$. Furthermore, Fig.~\ref{ex4sol} shows the numerical solutions of the primal and dual problems. In the primal problem, the convection velocity drives the solution diagonally from the bottom left to the top right of the domain. Consequently, a sharp boundary layer forms near the diagonal outflow boundary around $(1,1)$. A similar boundary layer appears in the dual problem, but with a reversed direction caused by the opposite convection velocity. The presented numerical solutions of the primal and dual problems are obtained on a mesh of $12288~(= 0.75 \cdot 128^2)$ quadrilateral elements. A comparison of the element-wise error contributions of $\eta_1^K$ and $\eta_2^K$ is shown in Fig.~\ref{ex4eta}, computed using a uniform quadrilateral mesh of $3072~(= 0.75 \cdot 64^2)$ elements. The results highlight that the error is most significant in the boundary layer region, where the solution exhibits rapid variations.


\section{Conclusion}

In this study, we have presented a VMS-based goal-oriented a posteriori error estimation framework for the OSGS stabilization method applied to CDR problems. This method was derived by selecting the SGS space orthogonal to the FE space -- a choice that has yielded reliable results. 

Two distinct approaches for a goal-oriented a posteriori error estimation have been presented. The first is an explicit VMS-based approach, in which the error is estimated by post-processing the FE solution. The second is an implicit duality-based approach that requires solving an additional adjoint problem to estimate the error. Numerical experiments confirm the efficient performance of the proposed error estimates and demonstrate that both methods produce similar local error estimates. However, the explicit approach is computationally more efficient than the implicit approach. We have also explained why goal-oriented error estimates based on global error norms are not appropriate for CDR problems.

Overall, the SGS-based error estimation approach has proven effective in accurately capturing numerical errors. This confirms that reliable error estimates can be modeled using VMS technology. In addition, the local error indicators identify higher error concentration regions in the domain, which motivates the development of mesh adaptation strategies. This represents a promising direction for future research work.

\section*{Acknowledgements}

R. Codina gratefully acknowledges the support received from the ICREA Acad\`emia Program, from the Catalan Government.
S. A. Khan acknowledges the financial support received from the Program Contract between CIMNE and the Catalan Government for the period 2021-2024.

\section*{References}

\bibliographystyle{elsarticle-num} 

\begin{thebibliography}{10}
\expandafter\ifx\csname url\endcsname\relax
  \def\url#1{\texttt{#1}}\fi
\expandafter\ifx\csname urlprefix\endcsname\relax\def\urlprefix{URL }\fi
\expandafter\ifx\csname href\endcsname\relax
  \def\href#1#2{#2} \def\path#1{#1}\fi

\bibitem{codina2000stabilization}
R.~Codina, Stabilization of incompressibility and convection through orthogonal sub-scales in finite element methods, Computer methods in applied mechanics and engineering 190~(13-14) (2000) 1579--1599.

\bibitem{hughes1995multiscale}
T.~J. Hughes, Multiscale phenomena: Green's functions, the {Dirichlet}-to-{Neumann} formulation, subgrid scale models, bubbles and the origins of stabilized methods, Computer methods in applied mechanics and engineering 127~(1-4) (1995) 387--401.

\bibitem{hughes1998variational}
T.~J. Hughes, G.~R. Feij{\'o}o, L.~Mazzei, J.-B. Quincy, The variational multiscale method—a paradigm for computational mechanics, Computer methods in applied mechanics and engineering 166~(1-2) (1998) 3--24.

\bibitem{hughes2018multiscale}
T.~J. Hughes, G.~Scovazzi, L.~P. Franca, Multiscale and stabilized methods, Encyclopedia of computational mechanics second edition (2018) 1--64.

\bibitem{codina2018variational}
R.~Codina, S.~Badia, J.~Baiges, J.~Principe, Variational multiscale methods in computational fluid dynamics, Encyclopedia of computational mechanics (2018) 1--28.

\bibitem{hauke2023review}
G.~Hauke, D.~Irisarri, A review of {VMS} a posteriori error estimation with emphasis in fluid mechanics, Computer Methods in Applied Mechanics and Engineering (2023) 116341.

\bibitem{ainsworth1997posteriori}
M.~Ainsworth, J.~T. Oden, A posteriori error estimation in finite element analysis, Computer methods in applied mechanics and engineering 142~(1-2) (1997) 1--88.

\bibitem{chamoin2023introductory}
L.~Chamoin, F.~Legoll, An introductory review on a posteriori error estimation in finite element computations, SIAM Review 65~(4) (2023) 963--1028.

\bibitem{gratsch2005posteriori}
T.~Gr{\"a}tsch, K.-J. Bathe, A posteriori error estimation techniques in practical finite element analysis, Computers \& structures 83~(4-5) (2005) 235--265.

\bibitem{hauke2006multiscale}
G.~Hauke, M.~H. Doweidar, M.~Miana, The multiscale approach to error estimation and adaptivity, Computer Methods in Applied Mechanics and Engineering 195~(13-16) (2006) 1573--1593.

\bibitem{hauke2006proper}
G.~Hauke, M.~H. Doweidar, M.~Miana, Proper intrinsic scales for a-posteriori multiscale error estimation, Computer methods in applied mechanics and engineering 195~(33-36) (2006) 3983--4001.

\bibitem{larson2005adaptive}
M.~G. Larson, A.~M{\aa}lqvist, Adaptive variational multiscale methods based on a posteriori error estimation: duality techniques for elliptic problems, Springer, 2005.

\bibitem{larson2007adaptive}
M.~G. Larson, A.~M{\aa}lqvist, Adaptive variational multiscale methods based on a posteriori error estimation: energy norm estimates for elliptic problems, Computer methods in applied mechanics and engineering 196~(21-24) (2007) 2313--2324.

\bibitem{bayona2018variational}
C.~Bayona-Roa, R.~Codina, J.~Baiges, Variational multiscale error estimators for the adaptive mesh refinement of compressible flow simulations, Computer Methods in Applied Mechanics and Engineering 337 (2018) 501--526.

\bibitem{colomes2018robustness}
O.~Colom{\'e}s, G.~Scovazzi, J.~Guilleminot, On the robustness of variational multiscale error estimators for the forward propagation of uncertainty, Computer Methods in Applied Mechanics and Engineering 342 (2018) 384--413.

\bibitem{hauke2008variational}
G.~Hauke, D.~Fuster, M.~H. Doweidar, Variational multiscale a-posteriori error estimation for multi-dimensional transport problems, Computer Methods in Applied Mechanics and Engineering 197~(33-40) (2008) 2701--2718.

\bibitem{codina2021posteriori}
R.~Codina, R.~Reyes, J.~Baiges, A posteriori error estimates in a finite element {VMS}-based reduced order model for the incompressible {Navier}-{Stokes} equations, Mechanics Research Communications 112 (2021) 103599.

\bibitem{baiges2017variational}
J.~Baiges, R.~Codina, Variational multiscale error estimators for solid mechanics adaptive simulations: an orthogonal subgrid scale approach, Computer Methods in Applied Mechanics and Engineering 325 (2017) 37--55.

\bibitem{babuvska1979analysis}
I.~Babu{\v{s}}ka, W.~C. Rheinboldt, Analysis of optimal finite-element meshes in {$R^1$}, Mathematics of computation 33~(146) (1979) 435--463.

\bibitem{kelly1983posteriori}
D.~W. Kelly, J.~De~SR~Gago, O.~C. Zienkiewicz, I.~Babuska, A posteriori error analysis and adaptive processes in the finite element method: {Part {I}}—error analysis, International journal for numerical methods in engineering 19~(11) (1983) 1593--1619.

\bibitem{babuvska1987feedback}
I.~Babu{\v{s}}ka, A.~Miller, A feedback finite element method with a posteriori error estimation: {Part {I}}. {The} finite element method and some basic properties of the a posteriori error estimator, Computer Methods in Applied Mechanics and Engineering 61~(1) (1987) 1--40.

\bibitem{babuvska1978posteriori}
I.~Babu{\v{s}}ka, W.~C. Rheinboldt, A-posteriori error estimates for the finite element method, International journal for numerical methods in engineering 12~(10) (1978) 1597--1615.

\bibitem{babuvvska1978error}
I.~Babu{\v{s}}ka, W.~C. Rheinboldt, Error estimates for adaptive finite element computations, SIAM Journal on Numerical Analysis 15~(4) (1978) 736--754.

\bibitem{jin1998posteriori}
H.~Jin, S.~Prudhomme, A posteriori error estimation of steady-state finite element solutions of the {Navier}—{Stokes} equations by a subdomain residual method, Computer methods in applied mechanics and engineering 159~(1-2) (1998) 19--48.

\bibitem{diez1998posteriori}
P.~D{\'\i}ez, J.~Egozcue, A.~Huerta, A posteriori error estimation for standard finite element analysis, Computer Methods in Applied Mechanics and Engineering 163~(1-4) (1998) 141--157.

\bibitem{larsson2010flux}
F.~Larsson, P.~D{\'\i}ez, A.~Huerta, A flux-free a posteriori error estimator for the incompressible {Stokes} problem using a mixed {FE} formulation, Computer methods in applied mechanics and engineering 199~(37-40) (2010) 2383--2402.

\bibitem{korotov-2008}
S.~Korotov, Global a posteriori error estimates for convection-reaction-diffusion problems, Applied Mathematical Modelling 32 (2008) 1579--1586.

\bibitem{zhang-et-al-2011}
Y.~Zhang, Y.~Hou, H.~Zuo, A posteriori error estimation and adaptive computation of conduction-convection problems, Applied Mathematical Modelling 35 (2011) 2336--2347.

\bibitem{john-novo-2013}
V.~John, J.~Novo, A robust {SUPG} norm a posteriori error estimator for stationary convection--diffusion equations, Computer Methods in Applied Mechanics and Engineering 255 (2013) 289--305.

\bibitem{du2021robust}
S.~Du, R.~Lin, Z.~Zhang, Robust recovery-type a posteriori error estimators for streamline upwind/{Petrov}--{Galerkin} discretizations for singularly perturbed problems, Applied Numerical Mathematics 168 (2021) 23--40.

\bibitem{sharma2021robust}
N.~Sharma, Robust a-posteriori error estimates for weak {Galerkin} method for the convection-diffusion problem, Applied Numerical Mathematics 170 (2021) 384--397.

\bibitem{tobiska2015robust}
L.~Tobiska, R.~Verf{\"u}rth, Robust a posteriori error estimates for stabilized finite element methods, IMA Journal of Numerical Analysis 35~(4) (2015) 1652--1671.

\bibitem{hauke2009variational}
G.~Hauke, D.~Fuster, Variational multiscale a posteriori error estimation for quantities of interest, Journal of Applied Mechanics 76 (2009) 021201.

\bibitem{granzow2017output}
B.~N. Granzow, M.~S. Shephard, A.~A. Oberai, Output-based error estimation and mesh adaptation for variational multiscale methods, Computer Methods in Applied Mechanics and Engineering 322 (2017) 441--459.

\bibitem{garg2019local}
V.~V. Garg, R.~H. Stogner, Local enhancement of functional evaluation and adjoint error estimation for variational multiscale formulations, Computer Methods in Applied Mechanics and Engineering 354 (2019) 119--142.

\bibitem{abdulle2013posteriori}
A.~Abdulle, A.~Nonnenmacher, A posteriori error estimates in quantities of interest for the finite element heterogeneous multiscale method, Numerical Methods for Partial Differential Equations 29~(5) (2013) 1629--1656.

\bibitem{wildey2008posteriori}
T.~Wildey, S.~Tavener, D.~Estep, A posteriori error estimation of approximate boundary fluxes, Communications in numerical methods in engineering 24~(6) (2008) 421--434.

\bibitem{li2022posteriori}
F.~Li, N.~Yi, A posteriori error estimates of goal-oriented adaptive finite element methods for nonlinear reaction--diffusion problems, Journal of Computational and Applied Mathematics 412 (2022) 114362.

\bibitem{kuzmin2010goal}
D.~Kuzmin, S.~Korotov, Goal-oriented a posteriori error estimates for transport problems, Mathematics and Computers in Simulation 80~(8) (2010) 1674--1683.

\bibitem{cnossen2006aspects}
J.~M. Cnossen, H.~Bijl, M.~I. Gerritsma, B.~Koren, Aspects of goal-oriented model-error estimation in convection-diffusion problems, in: ECCOMAS CFD 2006: Proceedings of the European Conference on Computational Fluid Dynamics, Egmond aan Zee, The Netherlands, September 5-8, 2006, Citeseer, 2006.

\bibitem{valseth2020goal}
E.~Valseth, A.~Romkes, Goal-oriented error estimation for the automatic variationally stable {FE} method for convection-dominated diffusion problems, Computers \& Mathematics with Applications 80~(12) (2020) 3027--3043.

\bibitem{chaudhry2014enhancing}
J.~H. Chaudhry, E.~C. Cyr, K.~Liu, T.~A. Manteuffel, L.~N. Olson, L.~Tang, Enhancing least-squares finite element methods through a quantity-of-interest, SIAM Journal on Numerical Analysis 52~(6) (2014) 3085--3105.

\bibitem{codina2000stabilized}
R.~Codina, On stabilized finite element methods for linear systems of convection--diffusion-reaction equations, Computer Methods in Applied Mechanics and Engineering 188~(1-3) (2000) 61--82.

\bibitem{codina1998comparison}
R.~Codina, Comparison of some finite element methods for solving the diffusion-convection-reaction equation, Computer methods in applied mechanics and engineering 156~(1-4) (1998) 185--210.

\bibitem{principe2010stabilization}
J.~Principe, R.~Codina, On the stabilization parameter in the subgrid scale approximation of scalar convection--diffusion--reaction equations on distorted meshes, Computer Methods in Applied Mechanics and Engineering 199~(21-22) (2010) 1386--1402.

\bibitem{codina2009subscales}
R.~Codina, J.~Principe, J.~Baiges, Subscales on the element boundaries in the variational two-scale finite element method, Computer Methods in Applied Mechanics and Engineering 198~(5-8) (2009) 838--852.

\bibitem{codina-gravenkamp-khan-2025}
R.~Codina, H.~Gravenkamp, S.~A. Khan, A posteriori analysis of the {OSGS} formulation for the convection-diffusion-reaction equation, (Submitted).

\bibitem{codina9}
R.~Codina, Stabilized finite element approximation of transient incompressible flows using orthogonal subscales, Computer Methods in Applied Mechanics and Engineering 191 (2002) 4295--4321.

\bibitem{codina2008analysis}
R.~Codina, Analysis of a stabilized finite element approximation of the {Oseen} equations using orthogonal subscales, Applied Numerical Mathematics 58~(3) (2008) 264--283.

\bibitem{Codina2001a}
R.~Codina, Pressure {{Stability}} in {{Fractional Step Finite Element Methods}} for {{Incompressible Flows}}, Journal of Computational Physics 170~(1) (2001) 112--140.
\newblock \href {https://doi.org/10.1006/jcph.2001.6725} {\path{doi:10.1006/jcph.2001.6725}}.

\bibitem{Castanar2020a}
I.~Casta{\~n}ar, J.~Baiges, R.~Codina, A stabilized mixed finite element approximation for incompressible finite strain solid dynamics using a total {{Lagrangian}} formulation, Computer Methods in Applied Mechanics and Engineering 368 (2020) 113164.
\newblock \href {https://doi.org/10.1016/j.cma.2020.113164} {\path{doi:10.1016/j.cma.2020.113164}}.

\bibitem{Gravenkamp2023a}
H.~Gravenkamp, R.~Codina, J.~Principe, A stabilized finite element method for modeling dispersed multiphase flows using orthogonal subgrid scales, Journal of Computational Physics 501 (2024) 112754.
\newblock \href {https://doi.org/10.1016/j.jcp.2024.112754} {\path{doi:10.1016/j.jcp.2024.112754}}.

\bibitem{Gravenkamp2023d}
H.~Gravenkamp, S.~Pfeil, R.~Codina, Stabilized finite elements for the solution of the {{Reynolds}} equation considering cavitation, Computer Methods in Applied Mechanics and Engineering 418 (2024) 116488.
\newblock \href {https://doi.org/10.1016/j.cma.2023.116488} {\path{doi:10.1016/j.cma.2023.116488}}.

\bibitem{diez-pares-huerta-2010}
P.~D\'{\i}ez, N.~Par\'es, A.~Huerta, Error estimation and quality control, Encyclopedia of Aerospace Engineering, {\rm Vol. 3, Part 15, Chapter 144} (2010) 1725--1734.

\bibitem{verfurth2005robust}
R.~Verf{\"u}rth, Robust a posteriori error estimates for stationary convection-diffusion equations, SIAM Journal on Numerical Analysis 43~(4) (2005) 1766--1782.

\bibitem{becker-rannacher-2001}
R.~Becker, R.~Rannacher, An optimal control approach to a posteriori error estimation in finite element methods, Acta Numerica 10 (2001) 1--102.

\bibitem{duprez-et-al-2020}
M.~Duprez, S.~P.~A. Bordas, M.~Bucki, H.~P. Bui, F.~Chouly, V.~Lleras, C.~Lobos, A.~Lozinski, P.-Y. Rohan, S.~Tomar, Quantifying discretization errors for soft tissue simulation in computer assisted surgery: A preliminary study, Applied Mathematical Modelling 77 (2020) 709--723.

\end{thebibliography}

\end{document}